\documentclass[12pt]{article}
\usepackage{amsfonts}
\usepackage{epsfig}
\title{Discrete Monodromy, Pentagrams,
and the Method of Condensation}
\author{Richard Evan Schwartz \thanks{\hskip 5 pt Supported by 
N.S.F. Research Grant DMS-0072607}}

\newtheorem{theorem}{Theorem}[section]

\newtheorem{lemma}[theorem]{Lemma}

\newtheorem{corollary}[theorem]{Corollary}

\def\startproof{{\bf {\medskip}{\noindent}Proof: }}

\def\endproof{$\spadesuit$  \newline}

\def\A{\mbox{\boldmath{$A$}}}% 
\def\C{\mbox{\boldmath{$C$}}}% 
\def\F{\mbox{\boldmath{$F$}}}% 
\def\N{\mbox{\boldmath{$N$}}}% 
\def\P{\mbox{\boldmath{$P$}}}% 
\def\R{\mbox{\boldmath{$R$}}}% 
\def\Z{\mbox{\boldmath{$Z$}}}% 

\begin{document}
\maketitle

\section{Introduction}

The purpose of this paper is to point out
some connections between:
\begin{enumerate}
\item The monodromy of periodic linear
differential equations;
\item The pentagram map, which
we studied in [{\bf S1\/}] and [{\bf S2\/}];
\item Dodgson's method of condensation for
computing determinants;
\end{enumerate}
We discovered most of these connections through
computer experimentation.

\subsection{Monodromy}

Consider the second order O.D.E.
\begin{equation}
\label{diff2}
f''(t)+\frac{1}{2} q(t) f(t)=0.
\end{equation}
Here $q(t)$ is $1$-periodic.
If $\{f_1,f_2\}$ is a basis for the solution space of
Equation \ref{diff2} then
there is some linear $T \in SL_2(\R)$ such
that $f_j(t+1)=T(f_j(t))$ for $j=1,2$.
The trace tr$(T)$, which is independent of
basis, is sometimes called the
{\it monodromy\/} of Equation
\ref{diff2}.
The ratio $f=f_1/f_2$ gives a smooth map
from $\R$ into the projective line.
Here $q$ is 
given by the {\it Schwarzian derivative\/}:
\begin{equation}
q=\frac{f'''}{f'}-\frac{3}{2} \bigg(
\frac{f''}{f'}\bigg)^2.
\end{equation}

Here is a discrete analogue of Equation \ref{diff2}.
The {\it cross ratio\/} of $4$ points
$a,b,c,d \in \R$ is given by
\begin{equation} \label{cro}
x(a,b,c,d)=\frac{(a-c)(b-d)}{(a-b)(c-d)}.
\end{equation}
A calculation shows that the quantity
\begin{equation}
\lim_{\epsilon \to 0} \frac{1}{\epsilon^2}
x(f(t-3 \epsilon),f(t-\epsilon),f(t+\epsilon),f(t+3 \epsilon))
\end{equation}
converges to a multiple of $q$, when $f$ is sufficiently smooth.
Thus, the cross ratio is a 
discrete analogue of the Schwarzian derivative.
Suppose  we have an infinite $n$-periodic sequence
$...q_n,q_1,q_2,...,q_n,q_1,...$.
We can find points $...,f_1,f_2,f_3,...$ in 
the projective line such that 
\begin{equation}
\label{diff1}
x(f_j,f_{j+1},f_{j+2},f_{j+3})=q_j \hskip 30 pt \forall j
\end{equation}
There will be a projective transformation
$T$ such that $f_{j+n}=T(f_j)$ for all $j$.
The conjugacy class of $T$ only depends on $q$.    To obtain
a numerical invariant, we can lift $T$
to $SL_2(\R)$ and take its trace.  
This quantity is a rational function in
the variables $q_1,...,q_n$. 

A main focus of this paper is a discrete analogue
for the third order case.  This analogue involves
infinite polygons in the projective plane.  In
analogy to the cross ratio 
we will define {\it projective invariants\/} of
polygons in \S 3.1.
We begin with an
infinite sequence $...,x_1,x_2,...$ of projective invariants
having period $2n$.  
These invariants determine, up to a projective
transformation, an infinite polygon which is
invariant under a projective transformation.
We call $P$ a {\it twisted $n$-gon\/}.
In other words, we have a map $P: \Z \to \R\P^2$
and a projective transformation $T$ such that
$P(n+j)=T(P(n))$ for all $j$.

The monodromies $\Omega_1$ and $\Omega_2$ corresponding to $T$
are rational functions of 
the variables $x_1,...,x_{2n}$. 
Let $[\cdot ]$ denote the floor function.
In \S 2.1 we will define polynomials
$O_1,...,O_{[n/2]},O_n$ and
$E_1,...,E_{[n/2]},E_n$.  We call these
polynomials the {\it pentagram invariants\/}.
We will express the monodromies explicitly
in terms of the pentagram invariants:
\begin{equation}
\label{main}
\Omega_1=\frac{(\sum_{k=0}^{[n/2]} O_k)^3}{O_n^2E_n};  \hskip 40 pt
\Omega_2=\frac{(\sum_{k=0}^{[n/2]} E_k)^3}{E_n^2O_n}.
\end{equation}

\subsection{The Pentagram Map}

Roughly, the {\it pentagram map\/} is the map which
takes the polygon $P$ to the polygon $P'$,
as indicated in Figre 1.  In \S 4 we
will give a precise definition, which
expresses the pentagram map as a
composition of two involutions
$\alpha_1$ and $\alpha_2$.

\begin{center}
\psfig{file=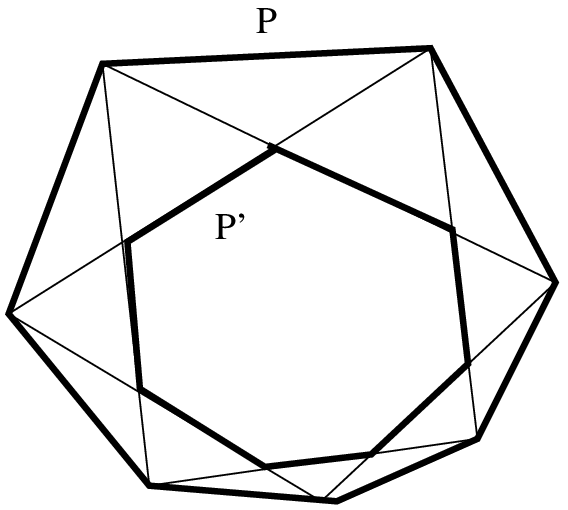}
Figure 1
\end{center}

Expressed in our
projective invariant coordinates$-$the cross
ratio generalizations discussed in the previous
section$-$the pentagram map
has the form
$\alpha_1(x_1,...,x_{2n})=(x'_1,...,x'_{2n})$ and
$\alpha_2(x_1,...,x_{2n})=(x''_1,...,x''_{2n})$ where
\begin{eqnarray}
\label{basic}
x_{2k-1}'=x_{2k} \frac{1-x_{2k+1}x_{2k+2}}{1-x_{2k-3}x_{2k-2}};
\hskip 25 pt
x_{2k}'=x_{2k-1}\frac{1-x_{2k-3}x_{2k-2}}{1-x_{2k+1}x_{2k+2}};
\cr \cr \cr
x_{2k+1}''=x_{2k}\frac
{1-x_{2k-2}x_{2k-1}}
{1-x_{2k+2}x_{2k+3}}
\hskip 25 pt
x_{2k}''=x_{2k+1}\frac
{1-x_{2k+2}x_{2k+3}}
{1-x_{2k-2}x_{2k-1}}
\end{eqnarray}
In these formulas, the indices are taken mod $2n$.
We let $\alpha=\alpha_1 \circ \alpha_2$.
In general, $\alpha$ has infinite order.

It turns out that the pentagram invariants are
invariant polynomials for the {\it pentagram map\/},
when it is expressed in suitable coordinates.

\begin{theorem} 
\label{trans}
$O_k \circ \alpha_j=E_k$ and
$E_k \circ \alpha_j=O_k$ for
$j=1,2$ and for all $k$.
\end{theorem}

\noindent
In \S 2 we will give a completely algebraic proof of
Theorem \ref{trans}.  In \S 3-4 we will give a
more conceptual proof which goes roughly as follows:
The pentagram map commutes with projective transformations
and therefore must preserve the monodromies
$\Omega_1$ and $\Omega_2$.
It follows from the general homogeneity
properties of Equation \ref{basic} that
the pentagram map must preserve the properly
weighted homogeneous pieces of the
monodromies, and these pieces are precisely
the pentagram invariants.   In \S 6 we prove

\begin{theorem}
\label{precise}
The pentagram invariants are algebraically
independent, so that
$\alpha$ has at least $2[n/2]+2$ algebraically independent
polynomial invariants.
\end{theorem}

We conjecture that the pentagram invariants give the complete
list of invariants for the pentagram map, at least when it
acts on the spaces of twisted $n$-gons.  We also
conjecturethat the algebraic varieties cut out
by the pentagram invariants are complex
tori, after a suitable compactification.
Finally we conjecture that the pentagram map acts
on these complex tori as a translation in
the natural flat metric. 

\subsection{The Method of Condensation}

Let $M$ be an $m \times m$ matrix. Let
$M_{NW}$ be the $(m-1) \times (m-1)$ minor
obtained by crossing off the last row and
column of $M$.  Here $N$ stands for ``north''
and $W$ stands for ``west''.  We define the
other three $(m-1) \times (m-1)$ minors $M_{SW}$,
$M_{NE}$ and $M_{SW}$ in the obvious way.
Finally, we define
$M_C$ to be the ``central''  $(m-2) \times (m-2)$ minor
obtained by crossing off all the extreme
rows and columns of $M$.  Dodgson's
identity says  \begin{equation}
\label{lc}
\det(M) \det(M_C)=\det(M_{NW}) \det(M_{SE})-
\det(M_{SW}) \det(M_{NE}).
\end{equation}
Assuming that $\det(M_C)$ is non-zero,
Equation \ref{lc} expresses $\det(M)$ 
as a rational function of determinants
of matrices of smaller size.  
This procedure can be iterated, expressing
the determinants of these smaller matrices
as rational functions of determinants of still
smaller matrices.  And so on.
This method of computing matrices
is called {\it Dodgson's method of
condensation\/}.  See
[{\bf RR\/}] for a detailed discussion of
this method and the rational functions that arise.

In \S 5 we will relate the pentagram map to
the method of condensation.   In some sense,
{\it the pentagram map
computes determinants\/}. We exploit this
point of view to prove

\begin{theorem}
\label{hyper}
Suppose that $P$ is a $4n$-gon whose sides
are alternately parallel to the $x$ and
$y$ axes.  Then (generically) the
$(2n-2)$nd iterate of the pentagram map transforms
$P$ into a polygon whose odd vertices are 
all collinear and whose even vertices are all
collinear.  
\end{theorem}

The surprise in Theorem \ref{hyper} is that $P$ could
have trillions of sides.
The pentagram map goes about its business for
trillions of iterations and then the whole thing
collapses all at once into a polygon whose
vertices lie on a pair of lines. 
Theorem \ref{hyper} is closely related
to the main result in [{\bf S3\/}],
which we proved by geometric methods.

\subsection{Paper Overview}
\begin{tabular}{ll} 
 {\bf \S 2: The Invariants\/}& \\
\S 2.1: Basic Definitions \\
\S 2.2: Proof of Theorem \ref{trans} \\ 

 {\bf \S 3: Discrete Monodromy\/}& \\
\S 3.1: PolyPoints and PolyLines &  \\
\S 3.2: Constructing the PolyPoints from its Invariants & \\
\S 3.3: The Final Calculation & \\

 {\bf \S 4: The Pentagram\/}& \\
\S 4.1: Basic Definitions& \\
\S 4.2: The Pentagram Map in Coordinates & \\
\S 4.3: Second Proof of Theorem \ref{trans} & \\
\S 4.4: Conic Sections \\

 {\bf \S 5: The Method of Condensation\/}& \\
\S 5.1: Octahedral Tilings & \\
\S 5.2: Picture of the Pentagram Map & \\
\S 5.3: Circulent Condensations & \\
\S 5.4: The Lifting Problem \\
\S 5.5: Degenerate Polygons \\
\S 5.6: Proof of Theorem \ref{hyper}\\

{\bf \S 6: Proof of Theorem \ref{precise}\/} \\
\S 6.1: Proof modulo the Vanishing Lemma \\
\S 6.2: Proof of the Vanishing Lemma
\end{tabular}

\subsection{Acknowledgements}

I would like to thank Peter Doyle, Bill Goldman,
Pat Hooper,
Francois Labourie, and John Millson for interesting
conversations related to this work. 
\newpage

\section{The Invariants}

\subsection{Basic Definitions}

All our definitions depend on a fixed
integer $n \geq 3$.  We will sometimes
suppress $n$ from our notation.
Let $Z=\{1,2,3,...,2n\}$.
We think of the elements of $Z$ as
being ordered cyclically, so that
$2n$ and $1$ are consecutive.
Also, in our notation all our indices
are taken cyclically.

We say that an {\it odd unit\/} of $Z$
is a subset having one of the two
forms:
\begin{enumerate}
\item $U=\{j\}$, where $j$ is odd.
\item $U=\{k-1,k,k+1\}$, where $k$ is even.
\end{enumerate}
We say that two odd units $U_1$ and $U_2$ are
{\it consecutive\/} if the set of odd
numbers in the union
$U_1 \cup U_2$ are consecutive.  For
instance $\{1\}$ and $\{3,4,5\}$ are
consecutive whereas
$\{1,2,3\}$ and $\{7,8,9\}$ are not.

We say that an {\it odd admissible subset\/} is
a nonempty subset $S \subset X$ consisting of a
finite union of odd units, no two of which
are consecutive.   We define the {\it weight\/}
of $S$ to be the number of odd units it contains.
We denote this quantity by $|S|$.   We define the
{\it sign\/} of $S$ to be the $+1$ is $S$ contains
an even number of singleton units, and $-1$ if $S$ contains
an odd number of singleton units.
As an example, the subset
$$\{1,5,6,7,11\}=\{1\} \cup \{5,6,7\} \cup \{11\}$$
is an odd admissible subset
if $n \geq 7$.  This subset has weight $3$ and sign $+1$.
As an exception to this rule, we call the
set $\{1,3,5,7,...,2n-1\}$ odd admissible as well.

Each odd admissible subset $S$ defines a
monomial $O_S \in R$:
\begin{equation}
O_S={\rm sign\/}(S) \prod_{j \in S} x_j.
\end{equation}
Let $O(k)$ denote the set of
weight $k$ odd admissible subsets of $Z$.
If $n$ is even then $O(k)$ is nonempty iff
$k \in \{1,2,...,n/2,n\}$.  If $n$ is odd then
$O(k)$ is nonempty iff $k \in \{1,2,...,(n-1)/2,n\}$.
We define
\begin{equation}
O_k=\sum_{S \in O(k)} O_S.
\end{equation}
By convention we set $O_0=1$.

We can make all the same definitions with the
word {\it even\/} replacing the word {\it odd\/}.
This leads to the definition of the $E$ polymonials.

\subsection{Proof of Theorem \ref{precise}}

Let $\alpha=\alpha_1 \circ \alpha_2$ be as
in the introduction.
For any rational function $f$, we define
$\alpha(f)=f \circ \alpha$.

By definition
\begin{equation}
O_n=x_1x_3...x_{2n-1}; \hskip 20 pt
E_n=x_2x_4...x_{2n}.
\end{equation}
If is easy to see directly from
Equation \ref{basic} that 
$\alpha_j(O_n)=E_n$ and
$\alpha_j(E_n)=O_n$.   When $n$ is even, we have
\begin{equation}
O_{n/2}=x_1x_5x_9...+x_3x_7x_{11}...; \hskip 20 pt
E_{n/2}=x_2x_6x_{10}...+x_4x_8x_{12}....
\end{equation}
Once again, it is easy to see directly from
Equation \ref{basic} that
$\alpha_j(O_{n/2})=E_{n/2}$ and
$\alpha_j(E_{n/2})=O_{n/2}$.
The interesting cases, which we now consider,
are when $k<n/2$.  We will show that
$\alpha_1(O_k)=E_k$.  The other cases have
similar derivations.

Before we treat the general case we consider an
example:  We have
$$O_1=\sum_{j=1}^n (-x_{2j+1}+x_{2j-1}x_{2j}x_{2j+1}).$$
Here indices are taken mod $2n$.
We compute easily that  \begin{equation}
\label{action}
\alpha_1(x_{2j+1}- x_{2j-1}x_{2j}x_{2j+1})=
x_{2j+2}-x_{2k+2}x_{2j+3}x_{2j+4}.
\end{equation}
Therefore
$$\alpha_1(O_1)=
\sum_{j=1}^n (-x_{2j+2}+x_{2k+2}x_{2j+3}x_{2j+4})=E_1.$$
This example suggests that the key to proving
Theorem \ref{precise} lies in partitioning
our polynomials in the right way.

Recall that $O(k)$ is the collection of weight $k$
odd admissible sequences.   Let $O_s(k) \subset O(k)$
consist of those sequences whose individual
units are singletons.  For instance
$\{1,5,9\}$ is a member of $O_s(3)$ as long as $n \geq 6$.
We call the odd units $\{j\}$ and $\{j-2,j-1,j\}$
{\it right partners\/}.   We say that a sequence
$S' \in O(k)$ is a {\it right partner\/} of
a sequence $S \in O_s(k)$ if every odd
unit in $S$ has a right partner odd unit in $S'$ and
{\it vice versa\/}.  For instance,
the sequence $S=\{1,7,19\}$ and $S'=\{1,5,6,7,17,18,19\}$
are right partners.   Also, $S$ is a right partner with
$\{1,5,6,7,19\}$.
For any $S \in O_s(k)$ let 
$S_R \subset O(k)$ be the subcollection of right
partners.   Every element of $O(k)$ has a unique
right partner in $O_s(K)$.  Therefore, we have a
partition
$$O(k)=\bigcup_{S \in O_s(k)} S_R.$$
Correspondingly, we can write
$$O_k=\sum_{S \in O_s(k)} RO_S; \hskip 30 pt
RO_S=\sum_{S' \in S_R} O_{S'}.$$

We can make all the same definitions, with
{\it even\/} replacing {\it odd\/} and
{\it left\/} replacing {\it right\/}.
Thus, we have a partition
$$E(k)=\bigcup_{S \in E_s(k)} S_L.$$
Here $S_L$ consists of the set of left
partners of $S$.
Correspondingly we can write
$$E_k=\sum_{S \in E_s(k)} LE_S; \hskip 30 pt
LE_S=\sum_{S' \in S_L} E_{S'}.$$

Below we will prove
\begin{lemma}
\label{seq}
For any sequence $S \in S_s(k)$ there is a sequence
$\overline S \in S_s(k)$ such that
$\alpha_1(RO_S)=LE_{\overline S}$. 
\end{lemma}

Since
$\alpha_1$ is an involution
the assignment $S \to \overline S$ is a bijection.
Summing over the individual terms we have
$\alpha(O_k)=E_k$.  Thus, Lemma \ref{seq} implies
Theorem \ref{precise}.  

To prove Lemma \ref{seq} we will decompose sequences
in $S_s(k)$, which could be quite complicated,
into much simpler sequences. We say that an
odd admissible sequence is {\it tight\/} if
it has the form
$$\{j,j+4,j+8,...,j+4a\}.$$  As usual, these numbers are
taken mod $2n$.  
Given any $S \in S_s(K)$ let $T_S$ denote the
set of maximal tight subsequences of $S$.
For instance, if $S=\{3,7,19,23,35\}$ and $n=18$ then
$T_S$ consists of the two sequences $\{19,23\}$
and $\{35,3,7\}$.  (The second sequence is
congruent mod $36$ to $\{35,39,43\}$.)
Since $k<n/2$ every sequence in $S_s(k)$ decomposes
in this way.

\begin{lemma}
\label{seq2}
If $S \in S_s(k)$ then
$RO_{S}=\prod_{T \in T_S} RO_T.$
\end{lemma}

\startproof
Let $T_S=\{T_1,...,T_h\}$.
Note that any number in $T_i$ is at least $5$ numbers
away from any number in $T_j$.  Otherwise,
$T_i \cup T_j$ would be tight, contradicting
maximality.  From this observation we see that
$T_1' \cup ... \cup T_h'$ is odd-admissible
for any choice of right partners
$T_1', T_2',...,T_h'$ of
$T_1,...,T_h$.  Conversely, any right
partner of $S$ decomposes this way.
Therefore
$S_R$ is precisely the union of the
sets of the form $T_1' \cup ... \cup T_h'$,
where $T'_j \in (T_j)_R$ is arbitrary.
Our lemma follows from this and from the
distributive law.
\endproof

\begin{lemma}
\label{seq3}
For any tight sequence $T$ there is a tight
sequence $\overline T$, having the same
length, such that
$\alpha_1(RO_T)=LE_{\overline T}$. 
\end{lemma}

\startproof
Let $P=RO_{T}$ and let $P'=\alpha(P)$.
Cyclically relabelling we can assume that
$T=\{3,7,11,...,4a+3\}$.
The set $T_R$ of right partners of $T$ consists of
$T$ and the sequence
$\{1,2,3,7,11,...,4a+3\}$.
Therefore
$$P=(1-x_1x_2)x_3x_7x_{11}...x_{4a+3}.$$
Writing $x_j'=\alpha(x_j)$ we have
$$P'=(1-x'_1x'_2)x'_3x'_7x'_{11}...x'_{4a+3}.$$
Using Equation \ref{action} we see that
$(1-x'_1x'_2)x'_3=x_4(1-x_5x_6)$.  We also have
$x_5x_6=x'_5x'_6$.
Therefore
$$P'=x_4(1-x'_5x'_6)x_7',x_{11}'...x'_{4a+3}.$$
Using Equation \ref{action} we see that
$(1-x'_5x'_6)x'_7=x_8(1-x_9x_{10})$.  We also have
$x_9x_{10}=x'_9x'_{10}.$
Therefore
$$P'=x_4x_8(1-x'_9x'_{10})x_{11}',...,x'_{4a+3}.$$
Continuing in this way we see that
$$P'=x_4x_8...x_{4a+4}(1-x'_{4a+5}x'_{4a+6})=
x_4x_8...x_{4a+4}(1-x_{4a+5}x_{4a+6}).$$
This last expression is
exactly $LE_{\overline T}$, where
$\overline T=\{4,8,...,4a+4\}$.
\endproof

Lemma \ref{seq} follows immediately from Lemma
\ref{seq2}, Lemma \ref{seq3}, and the uniqueness of
our decomposition.  As we mentioned before,
Theorem \ref{precise} follows from Lemma \ref{seq}.
This completes our proof.

\newpage

\section{Discrete Monodromy}

\subsection{PolyPoints and PolyLines}
\label{poly}

As in previous chapters we will fix some
positive integer $n \geq 3$.

Let $\P$ be the projective plane over the field
$\F$.   Say that a {\it PolyPoint\/} is a 
bi-infinite sequence $A=\{...A_{-3},A_1,A_5,..\}.$ of
points in $\P$.   (For technical reasons
we always index these points by integers
having the same odd congruence mod $4$.)
We assume also that there
is a projective transformation $T$ such
that $T(A_j)=A_{j+4n}$ for all $j \in \Z$.
We call $T$ the {\it monodromy\/} of $A$.

Say that a {\it PolyLine\/} is a 
bi-infinite sequence $B=\{...B_{-1},B_3,B_7,..\}$ of
lines in $\P$. 
We assume also that there
is a projective transformation $T$ such
that $T(B_j)=B_{j+4n}$ for all $j \in \Z$.
We call $T$ the {\it monodromy\/} of $B$.

Given two points $a,a' \in \P$ we let $(aa')$ be the
line containing these two points.  Given two
lines $b,b' \in \P$ we let $(bb')$ be the point
of intersection of these two lines.  Every PolyPoint
$A$ canonically determines a PolyLine $B$, by the
rule $B_j=(A_{j-2}A_{j+2})$.   At the same time
every PolyLine $B$ determines a PolyPoint $A$ by the
rule $A_j=(B_{j-2}B_{j+2})$.  In this case we call
$A$ and $B$ {\it associates\/}.   By construction
associates have the same monodromy.

The {\it dual space\/} to $\P$ is the space 
of lines in $\P$.  This space, denoted by
$\P^*$, is isomorphic to $\P$.
Indeed $\P^*$ is the projectivization of
the vector space dual to $\F^3$.
Any projective transformation $T: \P \to \P$
automatically induces a projective transformation
$T^*: \P^* \to \P^*$, and {\it vice versa\/}.
Any point in $\P$ canonically determines a
line in $\P^*$.  Likewise, points in $\P^*$
canonically determine lines in $\P$ and
lines in $\P^*$ canonically determine
points in $\P$.  The two spaces are on
an equal footing.

Given the PolyPoint $A$, we define
$A^*$ to be the PolyPoint in $\P^*$ whose
lines are given by the associate $B$.
If the points of $A$ are indexed by
numbers congruent to $1$ mod $4$ then
the points of $A^*$ are indexed by
numbers congruent to $3$ mod $4$, and
{\it vice versa\/}.
We make the same definitions for
PolyLines.
By construction
$A^{**}=A$ and $B^{**}=B$.  If $T$ is
the common monodromy of $A$ and $B$ then
$T^*$ is the common monodromy of
$A^*$ and $B^*$.  We call $A^*$ and
$B^*$ the {\it duals\/} of $A$ and $B$.

For any projective transformation $T$,
acting either on $\P$ or $\P^*$ we
define
\begin{equation}
\Omega_1(T)=\frac{{\rm tr\/}^3(\widetilde T)}{\det(\widetilde T)};
\hskip 30 pt
\Omega_2(T)=\Omega_1(T^*).
\end{equation}
Here $\widetilde T$ is a linear transformation whose
projectivization is $T$.  That is, $\widetilde T$ is a
{\it lift\/} of $T$.
It is easy to see that these quantities are
independent of lift.  Moreover,
$\Omega_j(T)$ only depends on the conjugacy
class of $T$.  Finally,
$\Omega_{3-j}(T^*)=\Omega_j(T)$ for any projective
transformation.

If $T$ is the monodromy of $A$ we call
$\Omega_1(T)$ and $\Omega_2(T)$ the
{\it monodromy invariants\/} of $A$.
By construction $A^*$ has the same
{\it set\/} of monodromy invariants as
$A$, but their order is switched.
The same goes for $B$.   If $S$ is
some other projective transformation,
then $A$ and $S(A)$ have the same
monodromy invariants.  Likewise,
$B$ and $S(B)$ have the same
monodromy invariants.

We now introduce our $2$-dimensional versions of the
cross ratio.   If $j$ is
one of the indices for the points of $A$ we define
\begin{eqnarray}
\label{invt}
p_{(j+1)/2}(A)=
x(A_{j + 8}, A_{j + 4}, (B_{j + 6} B_{j - 2}),
                            (B_{j + 6} B_{j - 6})) \cr \cr
q_{(j-1)/2}(A)=x(A_{j - 8}, A_{j - 4}, (B_{j - 6} B_{j + 2}),
                            (B_{j - 6} B_{j + 6}))
\end{eqnarray}
Here $x$ stands for the ordinary cross ratio,
as in Equation \ref{cro}.  
In the first equation, all $4$ points lie on the
line $B_{j+6}$.  In the second equation, all
$4$ points lie on $B_{j-6}$.   Conpare Figure 3 below.
If the points of $A$ are labelled by 
integers congruent to $1$ mod $4$ then the
invariants of $A$ are
$...q_0,p_1,q_2,p_3,...$.
If the points of $A$ are indexed by integers
congruent to $3$ mod $4$ then the invariants of
$A$ are
$...p_0,q_1,p_2,q_3,...$  In this chapter we will
only consider the case when the points of $A$ are
indexed by integers congruent to $1$ mod $4$,
though in the next chapter we will consider both
cases on an equal footing.

We can make all the same definitions for $B$, simply
by interchanging the two roles of $A$ and $B$ in
Equation \ref{invt}.
It turns out that our invariants are not
just invariant under projective transformations,
but also invariant under projective duality.
Precisely, we have
\begin{equation}
\label{dualinv}
p_j(A)=q_j(A^*); \hskip 15 pt
q_j(A)=p_j(A^*); \hskip 15 pt
p_j(B)=q_j(B^*); \hskip 15 pt
q_j(B)=p_j(B^*)
\end{equation}
for all relevant indices.  To see this symmetry,
we will consider an example.

Suppose that points of $A$ are labelled by
integers congruent to $1$ mod $4$.  The
first half of Figure 3 highlights the
$4$ points whose cross ratio is $p_3(A)$.
The second half shows the lines whose
cross ratio is used to define $q_3(A^*)$.
The highlighted
points are exactly the intersection points of
the highlighted line with an auxilliary line.
Hence, the two cross ratios are the same.

\begin{center}
\psfig{file=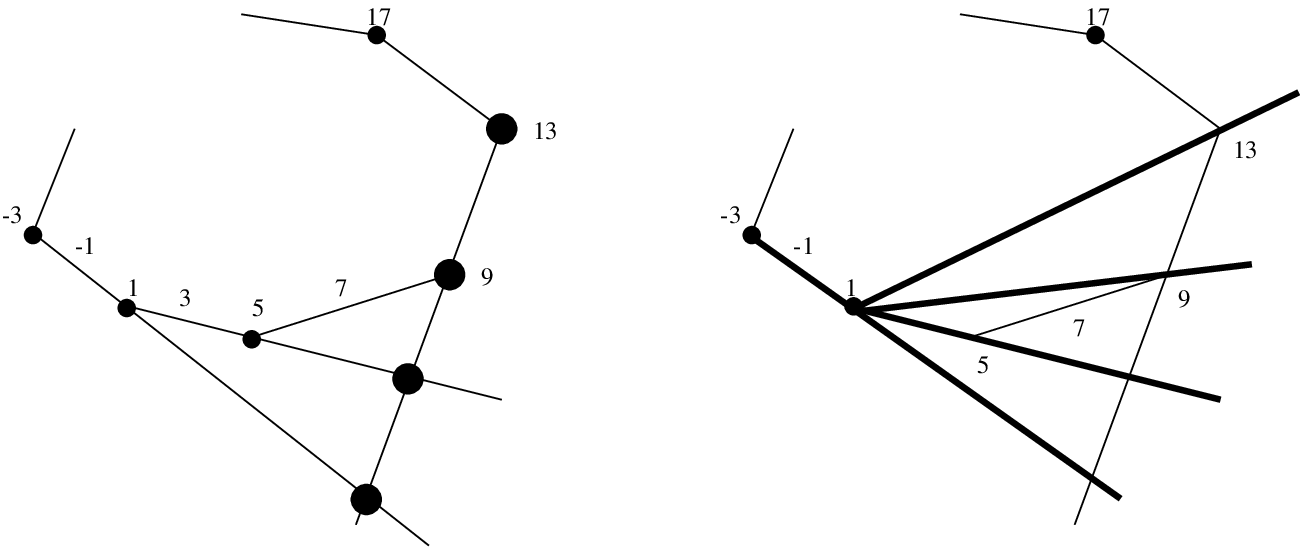}
Figure 3
\end{center}

This chapter is devoted to establising
Equation \ref{main}, which gives the formulas
for $\Omega_1$ and $\Omega_2$ in terms of
our invariants.  Given the formula for
$\Omega_1$, the formula for $\Omega_2$ follows
from projective duality and from 
Equation \ref{dualinv}.
Thus, to establish Equation \ref{main} it
suffices to derive the equation for
$\Omega_1$.

\subsection{Constructing the PolyPoint from its Invariants}
\label{monodromy series}

In \S 2 we constructed our polynomials from the variables
$x_1,...,x_{2n}$.  In this section we are going to
use the alternate list of variables
$p_1,q_2,p_3,q_4,...$.  The reason for the alternate
notation is that it is useful to distinguish the
even and odd variables in our constructions.
The polynomials in \S 2 are obtained from the ones
here using the substitution
$p_i \to x_i$ when $i$ is odd and
$q_i \to x_i$ when $i$ is even.

Suppose that $p_1,q_2,p_3,q_4,...$ are given variables.
We seek an infinite PolyPoint $A$ such that
\begin{equation}
\label{mainn}
p_{2i-1}(A)=p_{2i-1}; \hskip 15 pt
  q_{2i}(A)=q_{2i}; \hskip 15 pt i=1,2,3....
\end{equation}
What we mean by Equation \ref{mainn} is that we
wish to specify the points of $A$ in such a
way that the invariants we seek match a
specified list $p_1,q_2,p_3,...$.
Likewise, we seek a formula for the associate $B$.
For our purposes we only need the formulas
for ``half'' of $A$ and ``half'' of $B$.
That is, we just need to know $A_{-3},A_1,A_5,...$ and 
$B_{-5},B_{-1},B_3,...$.

Here we make the same definitions as in \S 2.1, with respect
$\Z$ (the integers) rather than the finite set $Z$.
To each admissible
sequence $S$ we associate a monomial $O_S$ in the formal
power series ring
$\A=\Z[[...p_1,q_2,p_3...]]$.
(Again, under the substitution mentioned above,
the ring $\A$ is identified with
$\Z[[...x_1,x_2,x_3,...]]$.) 
For instance if $S=\{1,2,3,9\}$ then
$O_S=-p_1q_2p_3q_9$.
We count the empty subset as both even and odd
admissible, and we define $O_{\emptyset}=E_{\emptyset}=1$.
Let $O$ be the
sum over all odd admissible sequences of finite
weight.  Likewise
let $E$ be the sum over all even admissible sequences
of finite weight. 
We have $O, E \subset \A$.
Given a pair of odd integers, $(r,s)$  we
define $O_r^s$  to be the polynomial
obtained from $O$ by setting $p_j$
equal to zero, for $j \leq r$ and $j \geq s$.
We make the same definitions, with
even replacing odd.

Let $A=\{...A_{-3}, A_{1}, A_{5},...\}$ and
$B=\{...B_{-5}, B_{-1}, B_3, ...\}$, where
(in homogeneous coordinates)

\begin{eqnarray}
\label{PA}
A_{-3}=[0,1,0]; \hskip 15 pt 
  A_{1}=[0,1,1]; \hskip 15 pt 
  A_{5}=[1,1,1];  \cr \cr
A_{4j+1}=[{O_1^{2j-1}},
\   {O_{-1}^{2j-1}+p_1 O_3^{2j-1}},
\   {O_{-1}^{2j-1}}]; \hskip 15 pt j=2,4,6... 
\end{eqnarray}

\begin{eqnarray}
\label{PB}
B_{-5}=[0,0,1]; \hskip 15 pt 
  B_{-1}=[1,0,0]; \hskip 15 pt 
  B_{3}=[0,1,-1]; \hskip 15 pt 
  B_{7}=[1,-1,0]; \cr \cr
B_{4j+3}=[{-E_2^{2j}+p_1q_2E_4^{2j}},
\   {E_0^{2j}},
\   {-E_0^{2j}+E_2^{2j}}]; \hskip 15 pt j=2,4,6...
\end{eqnarray}

In \S 5.2 we explicitly list out the first $7$ points of $A$.
We discovered these formulas as follows.  We normalized
the first few points of $A$ and then found the
equations for successive points using the definitions
of the invariants.  At some point we saw a pattern
in the growing polynomials we were generating.
The algebraic proofs we give in this section are
really more like verifications.  We did everything
on the computer and simply converted our observations
into a proof.  

The basic tool for us is the following set
of relations, which are easily derived.

\begin{eqnarray}
\label{relations}
O_r^s=0 \hskip 15 pt \forall r>s; \hskip 30 pt
  O_{s-2}^s=O_s^s=1; \cr
E_r^s=0 \hskip 15 pt \forall r>s; \hskip 30 pt
  E_{s-2}^s=E_s^s=1; \cr \cr
O_r^s=O_{r+2}^s
-p_{r+2}O_{r+4}^s+
P_{r+3} O_{r+6}^s; \hskip 15pt r <s. \cr \cr
E_r^s=E_{r+2}^s
-q_{r+2}E_{r+4}^s+
Q_{r+3}E_{r+6}^s. \hskip 15 pt r<s. \cr \cr
O_r^s=O_r^{s-2}-p_{s-2} O_r^{s-4}+P_{s-3} O_r^{s-6};
\hskip 15 pt r<s. \cr \cr
E_r^s=E_r^{s-2}-q_{s-2} E_r^{s-4}+Q_{s-3} E_r^{s-6};
\hskip 15 pt r<s. 
\end{eqnarray}
Here we have set \begin{equation} \label{shorten}
P_j=p_{j-1}q_jp_{j+1}; \hskip 30 pt Q_j=q_{j-1}p_jq_{j+1}.
\end{equation}

Let $\cdot$ stand for the dot product, and let
$\times$ stand for the cross product.

\begin{lemma}
\label{id2}
Let $k \geq 2$ and $d \geq 0$.  Then
$$A_{4k+1} \cdot B_{4k+3+4d}=
p_1q_2...q_{2k} E_{2k+2}^{2k+2d}.$$
$$B_{4k+3} \cdot A_{4k+5+4d}=
p_1q_2...q_{2k}p_{2k+1} O_{2k+3}^{2k+1+2d}.$$
\end{lemma}

\startproof
We will prove the first identity.  The second
one is very similar.
We use the notation 
$$(r,s)=O_{r}^{2k-1}E_s^{2k+2d}.$$
This notation should suggest to the reader that they plot the
various {\it points\/} $(r,s)$, given below, on a grid.
The result is a neat graphical representation of the
algebra we will be doing.
Using Equations \ref{PA} and \ref{PB} we have
\begin{equation}
\label{start}
A_{4k+1} \cdot B_{4k+3+4d}=
(-1,2)-(1,2)+p_1(3,0)+p_1q_2(1,4). \end{equation}
The basic relation
$$O_{-1}^*-O_1^*=-p_1 O_3^*+p_1q_2p_3O_5^*,$$
implies that
$$(-1,2) - (1,2) = -p_1(3,2)  + p_1q_2p_3(5,2).$$
Plugging this into Equation \ref{start} we have
\begin{equation}
\label{start1}
A_{4k+1} \cdot B_{4k+3+d}=p_1((3,0)-(3,2)+q_2(1,4)+q_2p_3(5,2).
\end{equation}
The basic relation
$$E_0^*-E_2^*=-q_2 E_4^*+q_2p_3q_4E_6^*,$$
implies that
$$(3,0)-(3,2)=-q_2(3,4)+q_2p_3q_4(5,4).$$
Plugging this into Equation \ref{start1} gives
\begin{equation} \label{start2}
A_{4k+1} \cdot B_{4k+3+4d}=
p_1q_2((1,4)-(3,4)+p_3(5,2)+p_3q_4(3,6)).
\end{equation}
note that Equation \ref{start2} has the same form
as Equation \ref{start} except that all the indices
have been shifted by $2$ and
a factor of $p_1q_2$ appears. 
This process repeats until we reach:
$$A_{4k+1} \cdot B_{4k+3+4d}=
p_1q_2...p_{2k-3}q_{2k-2} X,$$
where
$$X=(2k-3,2k+4)-(2k-1,2k+4)+
p_{2k-1}(2k+3,2k)+p_{2k-1}q_{2k}(2k-1,2k+2)$$ 
$$=1-1+0+p_{2k-1}q_{2k}=p_{2k-1}q_{2k}.$$
The last two equations combine to give
our identity.
\endproof

\begin{lemma}
\label{id1}
The following identities hold for all
$k \geq 2$.
\begin{enumerate}
\item $A_{4k+1} \times A_{4k+5}=
p_1...p_{2k-1} B_{4k+3}$. 
\item $B_{4k+3} \times B_{4k+7}=
q_2...q_{2k} A_{4k+5}$. 
\item $A_{4k+1} \cdot B_{4k+7}=p_1q_2...p_{2k-1}q_{2k}.$ 
\item $B_{4k+3} \cdot A_{4k+9}=p_1q_2...q_{2k}p_{2k+1}.$ 
\item $A_{4k+1} \cdot B_{4k+11}=p_1q_2...p_{2k-1}q_{2k}.$
\item $B_{4k+3} \cdot A_{4k+13}=p_1q_2...q_{2k}p_{2k+1}.$
\end{enumerate}
\end{lemma}

\startproof
We will derive these $6$ identities from the
two identities of Lemma \ref{id2}.
Taking $d=1$ we obtain Identities 3 and 4 listed above.
Taking $d=2$ we obtain Identities 5 and 6 listed above.
Taking $d=0$, we see that 
$$A_{4k+1} \cdot B_{4k+3}=
  B_{4k+3} \cdot A_{4k+5}=0.$$
Thus, we may write
$$A_{4k+1} \times A_{4k+5}=\lambda_k B_{4k+3},$$ for
some $\lambda_k$.
An easy calculation verifies that $\lambda_2=p_1p_3.$
Suppose, by induction, that $\lambda_{k-1}=p_1...p_{2k-3}$.
We use the fact that
$$(A_{4k+5} \times A_{4k+1}) \cdot A_{4k-3}=
  (A_{4k+1} \times A_{4k+3}) \cdot A_{4k-5},$$
and the already proven identities show that
$\lambda_{k}=p_1...p_{2k-1}$.  The case for the
$B$'s is similar. \endproof
\newline

Identity 1 
says that $A_i$ is the intersection
point of the lines $B_{i-2}$ and $B_{i+2}$.
Identity 2 says that
$B_{j}$ is the line determined by 
$A_{j-2}$ and $A_{j+2}$.   Hence $A$ and $B$ are
associates. 
Identities 3 and 4
say that generically $A_i \cdot B_{i-3} \not =0$ and
$B_j \cdot A_{j-3} \not =0$.
These statements also hold true
in the few cases where some of the indices are negative,
as may be verified by hand.  Thus, 
$A$ and $B$ are in general position for
generic choices of variables.

The first few identities of Equation $(*)$ can be
verified by hand.  We will show that
$p_{2k+5}(P)=p_{2k+5}$ for $k \geq 2$.  The case
for the $q$'s has a similar treatment.
To aid us in the computation, we recall some
vector identities:
\begin{eqnarray}
A \times (B \times C)=-(A \cdot B)C+(A \cdot C)B \cr \cr
  (A \times B) \times (A \times C)=((A\times B) \cdot C)A.
\end{eqnarray}
If $w_1,w_2,w_3,w_4$ are vectors lying in the same
$2$-dimensional linear subspace of $\F^3$ we define
\begin{equation}
\label{cr}
X(w_1,w_2,w_3,w_4)=\frac{(w_1 \times w_2)*(w_3 \times w_4)}
        {(w_1 \times w_3)*(w_2 \times w_4)}.
\end{equation}
The operation $*$ means coordinate-wise multiplication.
$X$ will be a vector of the form
$(x,x,x)$.   The number $x$ is the classical
cross ratio of the $4$ points in the projective
plane $\P$ represented
by the vectors $w_1,w_2,w_3,w_4$.

By definition
$$p_{2k+5}(P)= 
X(A_{4k+17},A_{4k+13},B_{4k+15} \times B_{4k+7},
		    B_{4k+15} \times B_{4k+3}),$$

First:
$$A_{4k+17} \times A_{4k+13}=
-p_1...p_{2k+5} B_{4k+15}.$$

Second:
$$(B_{4k+15} \times B_{4k+7}) \times 
  (B_{4k+15} \times B_{4k+3})=
  ((B_{4k+15} \times B_{4k+7}) \cdot B_{4k+3}) B_{4k+15}=$$
  $$((B_{4k+3} \times B_{4k+7}) \cdot B_{4k+15}) B_{4k+15}=
    (q_2...q_{2k} A_{4k+5} \cdot B_{4k+15}) B_{4k+15}=$$
  $$(q_2...q_{2k})(p_1,q_2,...,p_{2k+1}q_{2k+2}) B_{4k+15}.$$ 

Third:
$$A_{4k+17} \times (B_{4k+15} \times B_{4k+7})=$$
$$(A_{4k+17} \times B_{4k+7}) B_{4k+15}
 -(A_{4k+17} \times B_{4k+15}) B_{4k+7}=$$
$$-(p_1q_2...p_{2k+2}q_{2k+3}) B_{4k+15}.$$

Fourth:
$$A_{4k+13} \times (B_{4k+15} \times B_{4k+3})=$$
$$-(A_{4k+13} \times B_{4k+15}) B_{4k+3}+
  (A_{4k+13} \times B_{4k+7}) B_{4k+15}=$$
$$(p_1q_2...q_{2k}p_{2k+1}) B_{4k+15}.$$

Notice that all the terms are multiples of
the same vector.   
Using the formula for $X$, we have:
$$x=\frac{(p_1,p_{2k+5})(q_2...q_{2k})(p_1q_2,..p_{2k+1}q_{2k+2})}
	 {q_{2k+2} p_{2k+3}(p_1q_2,...,p_{2k+1})^2}=p_{2k+5}.$$

\subsection{The Final Calculation}
\label{monodromy}

Recall that $T$ is the monodromy of $A$.
We will first compute a lift
of $T$ to $GL_3(\F)$.
We may interpret an element in $GL_3(\F)$ as 
a triple $\widetilde T=(V_1,V_2,V_3)$ of vectors in $\F^3$.
The linear action of $T$ is then given as follows:
If $W=[w_1,w_2,w_3]$, then  \begin{equation}
\label{mon1}
\widetilde T(W)=w_1 V_1+w_2 V_2+ w_3 V_3.
\end{equation}

Consider the element
$\widetilde T=(V_1,V_2,V_3),$
where \begin{eqnarray}
\label{mon2}
V_1=p_1 A_{4n+5}-p_{2n+1} A_{4n+1}; \cr \cr
V_2=p_{2n-1}q_{2n}p_{2n+1} A_{4n-3}; \cr \cr
V_3= p_{2n+1}A_{4n+1}-
       p_{2n-1}q_{2n}p_{2n+1} A_{4n-3}. 
\end{eqnarray}

\begin{lemma}
$\widetilde T$ is a lift of $T$.
\end{lemma}

\startproof
Using Equations \ref{mon1}, \ref{mon2} and \ref{PA} we compute
$$\widetilde T(A_{-3})=\widetilde T[0,1,0]=
                     p_{2n-1}q_{2n}p_{2n+1} A_{4n-3}.$$
$$\widetilde T(A_{1})=\widetilde T[0,1,1]=
                     p_{2n+1} A_{4n+1}.$$
$$\widetilde T(A_{5})=\widetilde T[1,1,1]=
                     p_1 A_{4n+5}.$$
$$\widetilde T(A_9)=\widetilde T[1,1,1-p_1]=$$
$$\hskip 15 pt 
p_1(A_{4n+5}-p_{2n+1} A_{4n+1}+p_{2n-1}q_{2n}p_{2n-1} A_{4n-3})=
p_1 A_{4n+9}.$$
This last equality follows Equation \ref{relations},
applied componentwise to our vectors.
Our four computations show that the projectivization
of $\widetilde T$ has the same action on
the points $A_{-3}, A_{1}, A_{5}, A_9$
as $T$ does.  These $4$ points are (for generic
choice of variables) in
general position, and projective transformations
are determined by their action on $4$ general
position points. 
Hence the projective action of $\widetilde T$ coincides with
the action of $T$.
\endproof

Now we compute $\Omega_1$.  We
set 
$$\widetilde O=\sum_{i=0}^{[n/2]} O_k.$$
Before we make the next calculation we note
that $p_{2n+1}=p_1$, under the assumption that
the invariants are $2n$-periodic.   Now for
the calculation:
\begin{eqnarray}
\label{trace}
{\rm tr\/}(T)=V_{11}+V_{22}+V_{33}=  \cr \cr
(p_1 O_1^{2n+1}-p_{2n+1}O_1^{2n-1})+  \cr \cr
(p_{2n-1}q_{2n}p_{2n+1}O_{-1}^{2n-3}+  
  p_1p_{2n-1}q_{2n}p_{2n+1}O_{3}^{2n-3})+ \cr \cr
(p_{2n+1}O_{-1}^{2n-1}- 
p_{2n-1}q_{2n}p_{2n+1}O_{-1}^{2n-3})= \cr \cr
p_1\big([O_1^{2n+1}]+[O_{-1}^{2n-1}-O_1^{2n-1}]
      +[p_{2n-1}q_{2n}p_1 O_3^{2n-1}\big]) 
=p_1\widetilde O.
\end{eqnarray}

In the last line, we have bracketed terms
so as to isolate the different kinds of
terms in $\widetilde O$.   The first expression
describes the terms of $\widetilde O$ which involve $p_1$,
but not $p_{2n-1}q_{2n}p_1$.  The second expression
describes the terms of $\widetilde O$ which do not involve
$p_1$.  The third expression describes terms
of $\widetilde O$ which involve
$p_{2n-1}q_{2n}p_1.$

Again using the fact that $p_{2n+1}=p_1$ we have
\begin{eqnarray}
\label{dett}
\det(T)=(V_1 \times V_2) \cdot V_3= \cr \cr
(p_1)(p_{2n-1})(p_{2n-1}q_{2n}p_{2n+1})
((A_{4n-5} \times A_{4n-1}) \cdot A_{4n+3})= \cr \cr
p_1^3(p_3...p_{2n+1})^2(q_2...q_{2n})= 
p_1^3(p_1p_3...p_{2n-1})^2(q_2...q_{2n}).
\end{eqnarray}
Combining Equations \ref{trace} and
\ref{dett} we get the formula for
$\Omega_1$.  

\newpage

\section{The Pentagram}
\label{pentagram}

\subsection{Basic Definition}

We fix $n$, as in previous chapters. 
Let $P_1$ (respectively $P_2$) be the space of
PolyPoints which have $2n$ periodic
invariant coordinates, and whose points
are labelled by integers congruent to
$1$ (respectively $3$) mod $4$.  Likewise we define
$L_1$ and $L_2$ for PolyLines.

Define
\begin{equation}
X=P_1 \cup P_2 \cup L_1 \cup L_2.
\end{equation}
Suppose that $A=\{A_j\}$ is a PolyPoint.   We define
$\delta_1(A)=\{B_j\}$ where
\begin{equation}
B_j=(A_{j-2}A_{j+2}).
\end{equation}
$\delta_1(A)$ is just the associate of $A$.
We make the same definition for PolyLines.
$\delta_1$ is an involution of $X$ which
interchanges spaces $P_{1}$ and $L_{3}$ and 
interchanges 
$P_{3}$ and $L_{1}$.

  We define
$\delta_1(A)=\{B'_j\}$ where
\begin{equation}
B_j=(A_{j-4}A_{j+4}).
\end{equation}
$\delta_2$ is an involution of $X$ which
interchanges spaces $P_{1}$ and $L_{1}$ and 
interchanges 
$P_{3}$ and $L_{3}$.
Figure 5.1 shows the action of $\delta_2$ on a
PolyPoint in $P_1$.

\begin{center}
\psfig{file=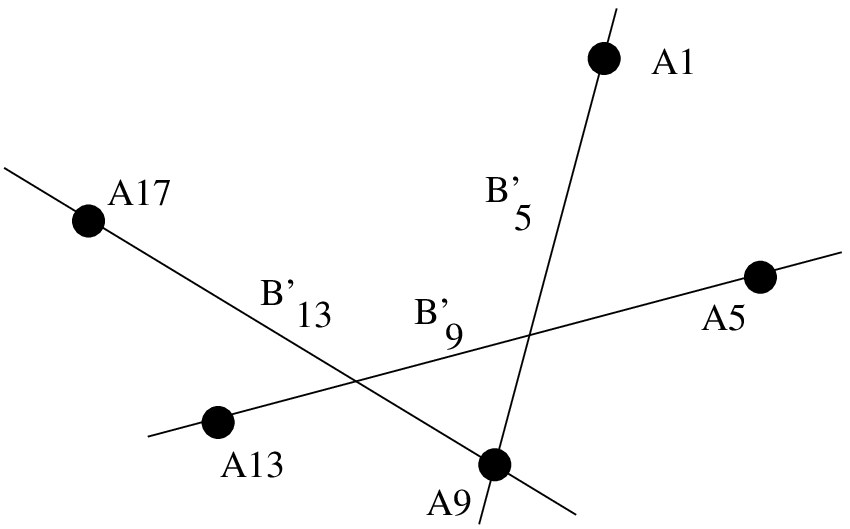}
Figure 5.1
\end{center}

We define
\begin{equation}
\alpha_1=\delta_1 \circ \delta_2 \circ \delta_1; \hskip 30 pt
\alpha_2=\delta_2.
\end{equation}

We call the pair $(\alpha_1,\alpha_2)$ the {\it pentagram map\/}.
Both $\alpha_1$ and $\alpha_2$ are involutions.  Moreover,
conjugation by $\delta_1$ interchanges these two maps.

\subsection{The Pentagram Map in Coordinates}
\label{coord}

There is a ``forgetful map'' which takes the variables
$p_j$ and $q_j$ and calls them $x_j$, regardless of
the letter.  There is also a map which takes a
polygon and sends it to its invariant coordinates.
Composing this coordinate map with the forgetful map
we get generically defined bijections from any
of our $4$ spaces to $\F^{2n}$.   In this way we
can work out the pentagram map in the same coordinates
used in Equation \ref{basic}.

We can exploit symmetry to reduce the amount
of computing we have to do.  First, it is
a consequence of Equation \ref{dualinv} and
the invariants of our coordinates under
projective duality that
\begin{enumerate}
\item The action of $\alpha_j$ on $P_j$ is the same
as the action of $\alpha_j$ on $L_j$, and this
action does not depend on $j$.
\item The action of $\alpha_j$ on $P_{3-j}$ is the same
as the action of $\alpha_j$ on $L_{3-j}$, and this
action does not depend on $j$.
\end{enumerate}
We will show that Equation \ref{basic} describes
the action of $\alpha_1$ and $\alpha_2$ on 
$P_1 \cup L_1$.  If we computed the action on
$P_2 \cup L_2$ we would have to interchange
$\alpha_1$ with $\alpha_2$ to get the
right formulas.  We will compute the action
of $\alpha_2$ on $P_1$.  The derivation for
$\alpha_1$ is similar and follows from symmetry.
To see that the formula in Equation \ref{basic} is correct,
we just have to work out a single invariant of
$\alpha_2(P)$.  The formulas for the other
invariants are forced by the dihedral
symmetry of our constructions.

We will not do this calculation by hand.  However, we will
set it up exactly, so that the inclined reader can type
everything into a symbolic manipulator and push the
button, as we did.  Let $A$ be as in
\S 4.  We will compute the invariant $q'_4=q_4(B')$, which
is the variable $x$ in the expression $X=(x,x,x)$, where
(as in Equation \ref{cr})
\begin{equation}
X=
X(B'_1, B'_5, A'_3 \times A'_{11},A'_3 \times A'_{15}).
\end{equation}
Here $A'$ is the associate of $B'$, so that
\begin{equation}
A'_3=B'_1  \times B'_5; \hskip 15 pt
A'_{11}=B'_{9} \times B'_{13}; \hskip 15 pt
A'_{15}=B'_{13} \times B'_{17}.
\end{equation}
Finally, we have
\begin{equation}
B'_1=A_{-3} \times A_5; \hskip 15 pt
B'_5=A_1 \times A_9; \hskip 15 pt ... \hskip 15 pt
B'_{17}=A_{13} \times A_{21}.
\end{equation}

To make our computation we need to know the 
$7$ points $A_{-3},...,A_{21}$. Once we have
these $7$ points, we just take a bunch of
cross products.

Using the convention of Equation \ref{shorten} we
list these points explicitly.
$$A_{-3}=[0,1,0]; \hskip 15 pt
A_{1}=[0,1,1]; \hskip 15 pt
A_{5}=[1,1,1]; \hskip 15 pt
A_9=[1,1,1-p_1];$$

$$A_{13}=[1-p_3,1-p_3+Q_2,1-p_1-p_3+Q_2];$$
$$A_{17}=[1-p_3-p_5+Q_4,1-p_3-p_5+Q_2+Q_4,
1-p_1-p_3-p_5+p_1p_5+Q_2+Q_4];$$

$$
A_{21}=[1-p_3-p_5-p_7+p_3p_7+Q_4+Q_6, $$
$$1-p_3-p_7-p_9+Q_2+Q_4+Q_6+p_3p_7 -p_7Q_2,$$
$$1-p_1-p_3-p_7+Q_2+Q_4+Q_6+p_1p_5+p_1p_7+p_3p_7
-p_1Q_6 - p_7 Q_2]$$
We double checked these formulas using
Equation \ref{relations}.  

When we plug
everything into Mathematica and compute,
we get
$$q'_4=p_5 \frac{1-q_6p_7}{1-q_2p_3}.$$
Forgetting about the lettering, we have
$$x'_4=x_5 \frac{1-x_6x_7}{1-x_2x_3}.$$
One can see that this exactly matches the
equation for $\alpha_2$ given in
Equation \ref{basic}, for $j=2$.
The general case follows from dihedral symmetry.

\subsection{Second Proof of Theorem \ref{precise}}

Suppose $A$ is a PolyPoint,
with invariants $p_1,q_2,...,q_{2n}$, as above.
Let $T$ be the projective transformation 
such that $T(A_{j+2n})=A_j$.  Let
$\Omega_1$ and $\Omega_2$ be the two
monodromy invariants of $A$.

We will just consider $\alpha_1$.  The proof
for $\alpha_2$ is the same.
Let $B'=\alpha_1(A)$.  Our constructions
commute with projective transformations.  Hence
$B'$ is also invariant under $T$.
This is to say that the dual of $B'$, which is
another PolyPoint, is
invariant under $T^*$, the dual of $T$.   Hence
$\alpha_1(\Omega_j)=\Omega_{3-j}$.
It now follows from
Equation \ref{main} that
\begin{equation}
 \alpha_1 \bigg(\frac{(\sum_{k=0}^{[n/2]} O_k)^3}{O_n^2E_n} \bigg)
=\frac{(\sum_{k=0}^{[n/2]} E_k)^3}{E_n^2O_n}.
\end{equation}
One can see easily from Equation \ref{basic} that
$\alpha_1(O_n)=E_n$ and
$\alpha_1(E_n)=O_n$.  Therefore
\begin{equation}
 \alpha_1\bigg(\sum_{k=0}^{[n/2]} O_k\bigg)=
\sum_{k=0}^{[n/2]} E_k.
\end{equation}

Now for the moment of truth.
Let $S_t: \F^{2n} \to \F^{2n}$ be as in
Equation \ref{hom}.  
Looking at
Equation \ref{basic} we see that
\begin{equation}
\alpha_1 \circ S_t=S_{t^{-1}} \circ \alpha_1.
\end{equation}
At the same time, we have
\begin{equation}
O_k \circ S_t=t^{-k} O_k;  \hskip 20 pt
E_k \circ S_t=t^{k} E_k.
\end{equation}
Therefore
\begin{equation}
 \alpha_1\bigg(\sum_{k=0}^{[n/2]} t^k O_k\bigg)=
\sum_{k=0}^{[n/2]} t^k E_k.
\end{equation}
Since this last equation is true for all $t$, we must
have $\alpha_1(O_k)=E_k$ for all $k$.
Since $\alpha_1$ is an involution,
$\alpha_1(E_k)=O_k$ for all $k$.  
This completes our proof.

\subsection{Conic Sections}

In this section we establish a technical result
used in later sections.  Namely, if $A$
is inscribed in a conic then $E_n(A)=O_n(A)$.
We continue using the notation established
above, and take $A \in P_1$.

\begin{lemma}
\label{conic1}
 Suppose $A$ is inscribed in a conic.
Then 
$$(1-q_{2k})(1-q_{2k-2}p_{2k-1})=
(1-p_{2k-1})(1-q_{2k}p_{2k+1});$$ $$
(1-p_{2k-1})q_{2k}(1-p_{2k+1})=
(1-q_{2k-2})p_{2k-1}(1-q_{2k})$$
for $k=1,...,n$.  Here indices are
taken mod $2n$.
\end{lemma}

\startproof
We will derive the first identity.  The second
one is a fairly straightforward rearrangement
of the first.   Figure 5.2 shows an application
of Pascal's theorem.  If the $6$ points
$A_{-3},A_1,A_5,A_9,A_{13},A_{17}$ lie on a
conic section then the points
$C_1,C_2,C_3$ lie on a line.

\begin{center}
\psfig{file=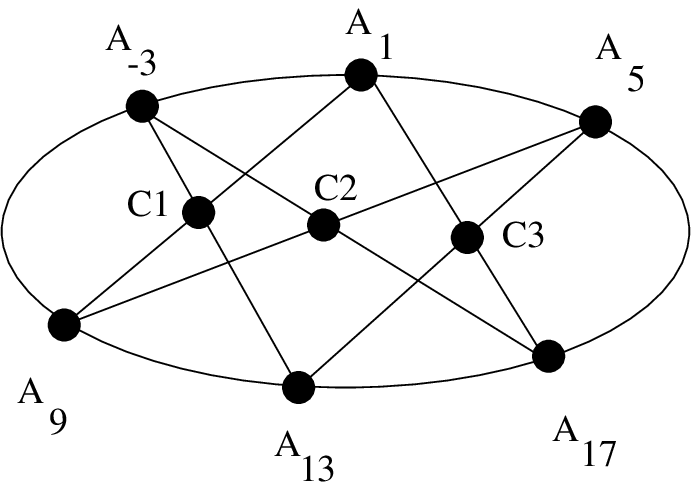}
Figure 5.2
\end{center}

If we express the $C_j$ in homogeneous coordinates,
we have \begin{eqnarray}
C_1=(A_{-3} \times A_{13}) \times (A_1 \times A_9) \cr \cr
C_2=(A_{-3} \times A_{17}) \times (A_5 \times A_9) \cr \cr
C_3=(A_{5} \times A_{13}) \times (A_1 \times A_{17})
\end{eqnarray}

The condition that the $C$s lie on the same line is
given by
\begin{equation}
\label{det}
\det(C_1,C_2,C_3)=0.
\end{equation}
In other words, we arrange these vectors into
a $3 \times 3$ matrix and set the determinant
equal to $0$.  

When we do this calculation on Mathematica, we see that
Equation \ref{det} holds if and only if
$(1-q_4)(1-q_2p_3)=(1-p_3)(1-q_4p_5).$  If all points
of $A$ lie on the same conic, then by symmetry the
same relation holds with the indices shifted
by $2,4,6...$
This completes our proof.
\endproof

\begin{corollary}
\label{later}
If $A$ is inscribed in a conic then
$O_n(A)=E_n(A)$.
\end{corollary}

\startproof
Taking the product of the first equation over all
$k$, and cancelling terms, we see that
$\prod (1-p_{2k-1})=\prod (1-q_{2j})$.
Taking the product of the second equation,
over all $k$, we get
$O_n (\prod (1-q_{2k})^2=
E_n (\prod (1-p_{2k-1})^2$.
Using the first equation to cancel terms,
we are left with $O_n=E_n$.
\endproof

\newpage

\section{The Method of Condensation}

\subsection{Octahedral Tilings}

There is a tiling $T$ of $\R^3$ by octahedra, which
can be described as follows.  The vertices of
$T$ have the form $(a,b,c) \in \Z^3$, subject
to the constraints that the three coordinates
are either all even or all odd.
Two vertices are joined by an edge if their
distance is exactly $\sqrt 3$.  One of the
octahedra in $T$ has the $6$ vertices
$(0,0,0)$,
$(\pm 1, \pm 1,1)$ and $(2,0,0)$.
All the other octahedra are translates of this one.
In our model octahedron, we call $(2,0,0)$ the
{\it top\/} and call $(0,0,0)$ the {\it bottom\/}.
We call $(-1,1,1)$ the {\it northwest vertex\/},
$(1,1,1)$ the {\it northeast vertex\/} and so forth.
We extend this definition to all octahedra
using translations.

Suppose that $M$ is an $m \times m$ matrix.
Dodgson's method of condensation involves the
connected square minors of $M$.  
Suppose that
\begin{equation}
H=\{M_{ij}|\ a_1 \leq i \leq a_2;\ b_1 \leq j \leq b_2\}
\end{equation}
is such a minor.  We define
$f(H)=(a,b,c)$ where
\begin{equation}
a=a_1+a_2; \hskip 20 pt b=b_1+b_2; \hskip 20 pt
c=a_1-a_2=b_1-b_2.
\end{equation}
For instance, if $H$ is the singleton $M_{11}$ then
$f(H)=(2,2,0)$.   If $H=M$ then
$f(H)=(m+1,m+1,m-1)$.  We let $[M]$ denote the
set of vertices of the tiling which have the
form $f(H)$.  Note that $[M]$ looks like
a pyramid with square base.

We can label the point $f(H) \in [M]$ by
$\det(H)$.   Dodgson's identity gives a
single relation for each octahedron:
\begin{equation}
V_t V_b= V_{nw} V_{se} - V_{sw} V_{ne}.
\end{equation}
Here $V_t$ is the label of the top vertex,
$V_b$ is the label of the bottom vertex, and
so forth.
To compute $\det(M)$, begins at the base of the
pyramid and computes successive layers using
the octahedron rule.  At the end, one arrives
at the apex of the pyramid, with the final answer.

What we do next works most gracefully when
$\det(M) \not =0$.  We say that a labelling
of a horizontal layer of $T$ is {\it constant\/}
if every vertex gets the same label.  
We say that a {\it sandwich
condensation\/} is a labelling of all the layers
of $T$ between two constant horizontal layers.
Here is an example:
We label the layer $\{z=0\}$ by $1$'s.
Next, we label the layer
$\{z=1\}$ so as to be doubly periodic, with period
$m$ in each direction.  In one $m \times m$ block, we put the
entries of $M$.
The second layer will be labelled
by determinants of $2 \times 2$ connected minors
of $M$ and its cyclic permutations.
By Dodgson's identity, the third layer
will be labelled by the determinants of
the $3 \times 3$ connected
minors of $M$ and its cyclic permutations.
And so forth.  The $m$th layer will be
labelled by the determinants of $M$ and
its cyclic permutations.  But all
these determinants have the same value.  This
is our sandwich.  It has {\it width\/} $m$.

To consider a more generic situation,
we say that a {\it condensation\/}
is a labelling of the vertices of the
tiling $T$, such that the labelling satisfies
the octahedron rule at every octahedron.
The labellings on two successive horizontal layers
determine the condensation.  If these labels
are chosen generically, in $\C$ for example,
then we will not
encounter singularities when trying to
propagate the condensation to the other
layers, both above and below the two
initial ones.

\subsection{Picture of the Pentagram Map}

\begin{center}
\psfig{file=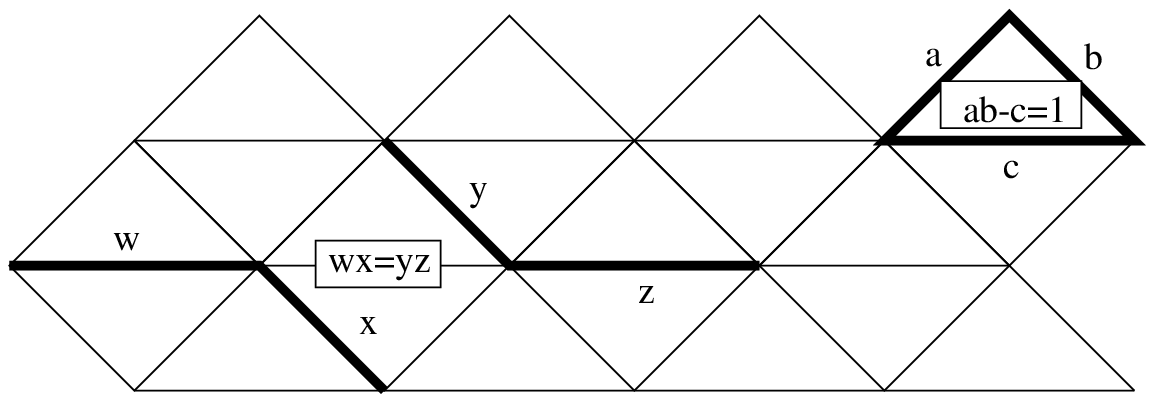}
Figure 5.1
\end{center}

Figure 5.1 shows a tiling of the
plane by isosceles triangles.  Suppose one labels
the edges of the tiling with elements of a field,
in a way which is invariant under a horizontal
translation, subject to the compatibility rules
indicated in the figure.    
(These compatibility rules are supposed
to hold for all configurations isometric to the
ones highlighted.)  We call such a labelling 
a {\it pentagram labelling\/}.

Given a pentagram labelling, the maps
$\alpha_1$ and $\alpha_2$,
the two involutions from Equation \ref{basic},
express how one one deduces the labellings on a given
row from the labellings on the rows above or
below it.  Thus, an orbit of the pentagram map
is encoded by a pentagram labelling.

\subsection{Circulent Condensations} 

Consider the linear projection
$\pi: \R^3 \to \R^2$ given by
\begin{equation} \label{proj}
\pi(x,y,z)=(2x - y,z).
\end{equation}
We say that a condensation is {\it circulent\/}
to $\pi$ if every vertex in a fiber
$\pi^{-1}(p)$ gets the same label.
Referring again to Dodgson's method of
condensation, the circulent labellings
correspond to certain circulent matrices.

\begin{center}
\psfig{file=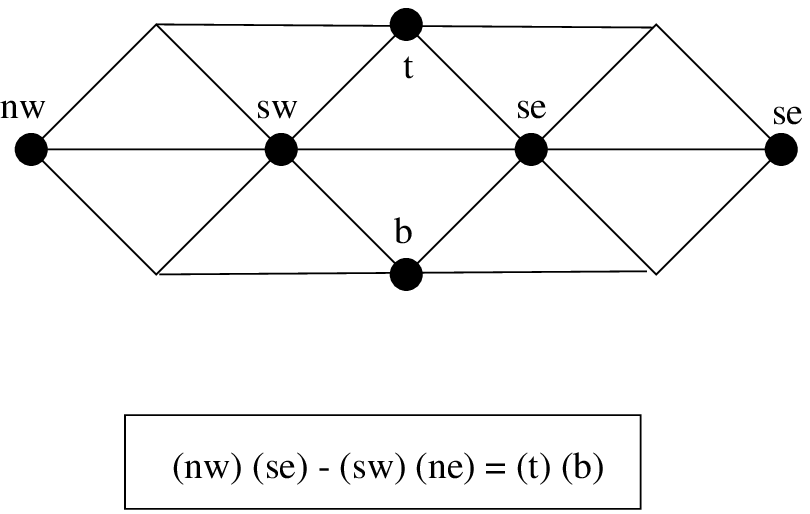}
Figure 5.2
\end{center}

We can use $\pi$ to push a circulent condensation
into the plane, without losing any information.
$\pi$ maps the vertices of $T$ to the vertices
of the tiling $\tau$ shown in Figure 5.1.
Figure 5.2 shows the image of a single octahedron under $\pi$,
superimposed onto $\tau$.
Figure 5.2 also shows the local rule satisfied
by the image of an circulent condensation.

We can think of an circulent condensation as a
labelling of the vertices of $\tau$ which
satisfies the local condition shown in Figure 5.2.
Figure 5.3 shows how to convert from an
circulent condensation to a pentagram labelling.

\begin{center}
\psfig{file=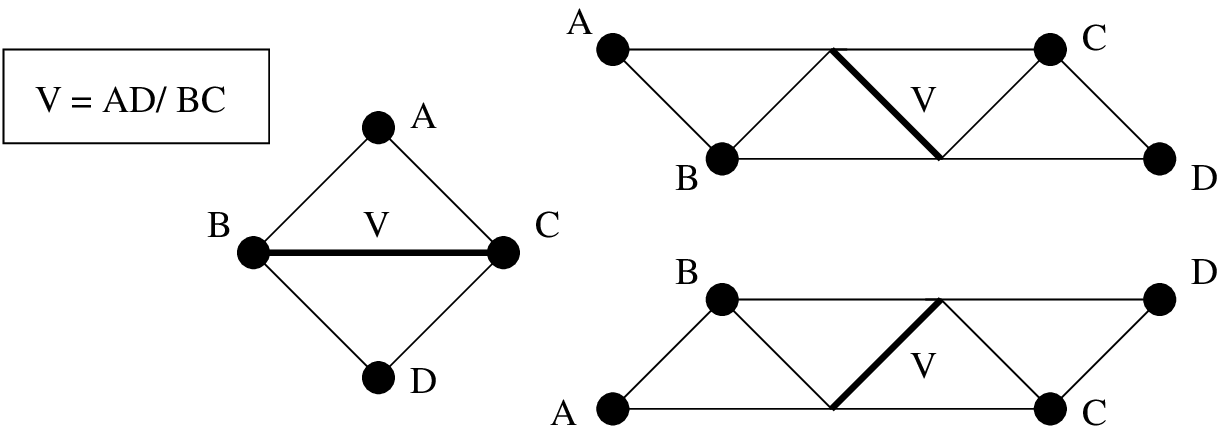}
Figure 5.3
\end{center}

To verify the first of the two compatibility
rules for pentagram labellings we refer to
Figure 5.4 and compute
$$ab-c=\frac{CD}{AF} \frac{AG}{CE} - \frac{BH}{EF}=
\frac{DG-BH}{EF}=1.$$
The last equality is the compatibility rule for the
circulent condensation.

\begin{center}
\psfig{file=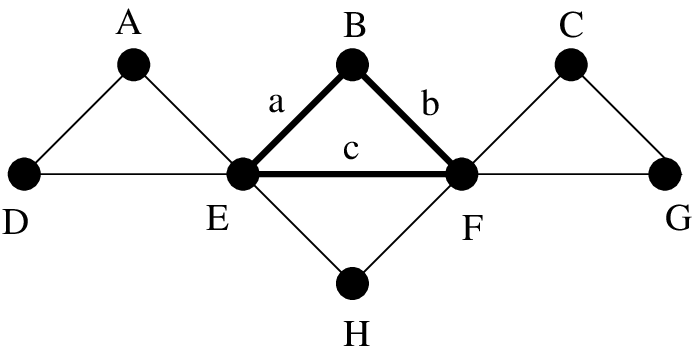}
Figure 5.4
\end{center}

To verify the second compatibility rule we refer to
Figure 5.5 and compute
$$wx=\frac{AI}{EF} \frac{EK}{GI}=\frac{AK}{FG}=
\frac{AH}{CF} \frac{CK}{GH}=yz.$$

\begin{center}
\psfig{file=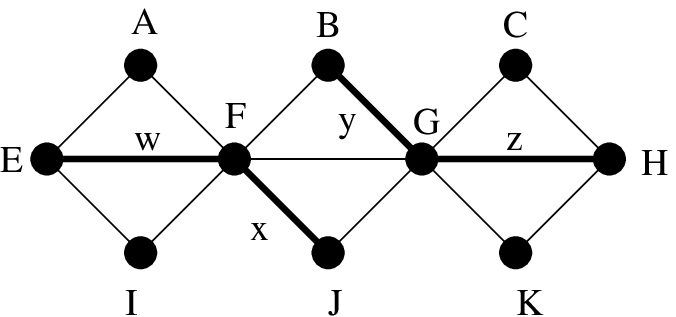}
Figure 5.5
\end{center}

\subsection{The Lifting Problem}

As usual, suppose that $n$ is fixed, as well as
a suitable base field.
Let $C$ denote the space of circulent condensations
which are periodic under a horizontal translation
which shifts the labels by $n$.  Let $P$ denote the
space of pentagram labels which are periodic under
the same translation.  The translation given in
the previous section gives a map $\psi: C \to P$.
In this section we ask about the extent to which
$\psi$ is invertible.

An element of $C$ is determined by the values it
attains on two successive horizontal rows.  
An element of $P$ is determined by the values
it takes on the non-horizontal edges of a single
row.  In both cases one puts down $2n$ values
before the period repeats.   That is, both
$C$ and $P$ are $2n$ dimensional.  The map
$\psi$ is a quotient map and certainly not
$1$-to-$1$.  For instance, all the scalar
multiples of a single element of $C$ get
mapped to the same element of $P$.
Since $C$ and $P$ have the same dimension and
$\psi$ is far from injective, $\psi$ is
also far from surjective.

To give a simple and clean answer we will assume
that $n=4m$ for some $m \in \N$.  The other cases
have a somewhat messier analysis.  Given an
element $X=(x_1,...,x_{2n}) \in P$ we have the
pentagram invariants defined in \S 2.1.

\begin{theorem}
\label{lift}
If $n=4m$ than an element
$X \in P$ is contained in the image of
$\psi$ if and only if
$O_{2m}(X)=E_{2m}(X)=2$ and $O_{n}(X)=E_{n}(X)=1$.
\end{theorem}

\startproof
For $j=1,2,3,4$ we define
\begin{equation}
\label{mod4}
f_j=\prod_{i \equiv j {\rm mod \/} 4} x_i.
\end{equation}
It is not hard to see that
$O_{2m}=f_1+f_3$.  Also $O_{n}=f_1f_3$.
By hypotheses, we have $f_1+f_3=2$ and
$f_1f_3=1$.  This forces $f_1=f_3=1$.
Likewise $f_2=f_4=1$.
Conversely, if $f_j=1$ for $j=1,2,3,4$ then
we have the hypothesis of this theorem.
All in all, the hypothesis of this theorem
are equivalent to the hypothesis that
$f_j=1$ for $j=1,2,3,4$.  We will work
with this latter hypothesis.

First we prove necessity.  Suppose
$X$ is in the image of $\psi$.
We focus on a single row of our triangulation.
Let $c_k$ be the label of the vertex which
is incident to the edges labelled 
$x_{k-1/2}$ and $x_{k+1/2}$.   Here $k$ is a half-integer.
We have
\begin{equation}
x_1=\frac{c_{-3/2} c_{7/2}}{c_{-1/2}c_{5/2}}; \hskip 15 pt
x_5=\frac{c_{5/2} c_{15/2}}{c_{7/2}c_{13/2}}; \hskip 15 pt
x_9=\frac{c_{13/2} c_{23/2}}{c_{15/2}c_{21/2}} \hskip 15 pt ...
\end{equation}
When we compute $f_1$, each $c$-term appears exactly
once in the numerator and exactly once in the
denominator.
Hence $f_1=1$.  A similar argument shows that
$f_2=f_3=f_4=1$.

To prove sufficiency, suppose that
$f_1=f_2=f_3=f_4=1$. We claim that there exist
labels $r_1,r_2,...,r_{2n}$ such that
\begin{equation}
x_j=\frac{r_{j+2}}{r_{j-2}}.
\end{equation}
To see this, we pick $r_1$ arbitrarily.
Then set, successively
\begin{equation}
r_5=x_3r_1; \hskip 15 pt
r_9=x_7x_3r_1; \hskip 15 pt
r_{11}=x_{11}x_7x_3 r_1 \hskip 15 pt...
\end{equation}
Since $n \equiv 0$ mod $4$ when we cycle through
one period we return back to the value
$r_1 f_3=r_1$. 
Thus, we can make a completely consistent
choice of $r$ for all indices congruent to
$1$ mod $4$.  The same argument works for the
other congruences.

Given the $r$'s, it suffices to find $c$'s such that
\begin{equation}
r_j=\frac{c_{j-1/2}}{c_{j+1/2}}.
\end{equation}
In our construction of the $r$'s we had a free
choice for each congruence mod $4$.
Thus, (in the generic case)
we can choose the $r$'s so that
\begin{equation}
\rho_j=\prod_{i \equiv j {\rm mod\/} 4} r_i=1
\end{equation}
for $j=1,2,3,4$.
Finding the $c$'s in terms of the $r$'s is exactly
the same problem as finding the $r$'s in terms of
the $x$'s.  The $\rho$'s play the same role
as the $f$'s.  Thus, we can find our $c$'s, which
gives a preimage of $X$ in $C$.
\endproof

\subsection{Degenerate polygons}

We fix $n=4m$.   To avoid trivial cases, we assume
that $m$ is large$-$say, $m \geq 3$.  We work
over $\C$.  Everything we do is understood to
be defined for the generic example, but
perhaps not for every example.  
The reader should insert this caveat before
every assertion.
We will work with PolyPoints
and PolyPoints which are periodic mod $n$.
By this we mean that the points or lines
themselves are periodic, not just the
projective invariants.  We will call such
objects {\it periodic PolyPoints\/} or
{\it periodic PolyLines\/}.  We take these
to be elements of the spaces $P_1$ and $L_1$
introduced in \S 5.1.

We say that the periodic PolyPoint $A=\{...A_1,A_5,...\}$ is
{\it degenerate\/} if $...A_1,A_9,A_{13}...$ lie
on the same line and if $...A_3,A_{11},A_{15},...$ all
lie on the same line.  We make the dual definition
for PolyLines.  A polygon (in the ordinary sense)
satisfying the hypothesis of Theorem \ref{hyper}
naturally determines a degenerate PolyLine.
One simply considers the lines considered by the
edges of the polygon.

\begin{lemma}
\label{recip}
A periodic element in $P_1$ is degenerate iff
its projective invariants satisfy $q_jp_{j+1}=1$
for all $j$.  Dually, a periodic element in $L_1$
is degenerate if and only if its projective
invariants satisfy $p_jq_{j+1}=1$ for all $j$.
\end{lemma}

\startproof
By duality it suffices to prove this result for
PolyPoints.  Referring to the PolyPoint constructed
in \S 4, we compute
$$\det(A_{-1},A_5,A_{13})=p_1(1-q_2p_3).$$
The three points in question lie on the same
line iff $q_2p_3=1$.  By symmetry, similar
statements hold for other triples of points.
\endproof

\begin{lemma}
\label{closed}
If $A$ is a periodic $2n$-gon then
$$(\sum_{k=0}^{[n/2]} O_k)^3=27 O_n^2E_n; \hskip 15 pt
(\sum_{k=0}^{[n/2]} E_k)^3=27 E_n^2O_n.$$
\end{lemma}

\startproof
This is a corollary of Equation \ref{main}. If
$A$ is $2n$-periodic then its monodromy
$T$ is the identity.  Likewise,
$T^*$ is the identity matrix.
Hence $\Omega_1=\Omega_2=27$.
\endproof

For Theorem \ref{hyper} we only
need the second statement of the next result.
The spaces $P_1$ and $L_1$ are defined in
\S 5.1.

\begin{lemma}
\label{van}
A degenerate element in either $P_1$ or $L_1$
has pentagram invariants \begin{enumerate}
\item $O_k=E_k=0$ for $k<2m$ and 
\item $O_{2m}=E_{2m}=2$ and
$O_n=E_n=1$.
\end{enumerate}
\end{lemma}

\startproof
First suppose that $k<2m=n/2$.  Using
the reasoning in the proof of Theorem \ref{trans} we
can group
$O_k$ into pieces having the form
$\prod \mu$, where each $\mu$ is a monomial
times  $(1-q_jp_{j+1})$ for some $j$.
(Compare Lemma \ref{seq2}, but interchange
$p$ with $q$.)
By Lemma \ref{recip} each $\mu$ is zero
and hence $O_k=0$.  The same argument
works for $E_k$.

For the second item note that
$O_nE_n=(p_2q_3)(p_4q_5)...=(1)(1)...=1.$
Elements in $P_1$ or $L_1$ are degenerate
limits of polygons which are either inscribed
in, or superscribed about, a conic section.
Thus, Lemma \ref{later} applies to these elements,
by continuity.
It now follows from Lemma \ref{later}
that $E_n^2=O_n^2=1$.   If $O_n=-1$, for any particular
example, then by continuity
$O_n$ is identically $-1$.  However, we
can choose a $4n$-point, contained in $\R^2$, which is essentially
the double cover of a $2n$-point.  The
invariants of such a $4n$-point would
repeat with period $2n$, forcing $O_n$ to be
a square of a real number.  This rules
out the possibility that $O_n=-1$.  Hence
$O_n=1$ for every example.  Likewise
$E_n=1$ for every example.
If now follows from Corollary \ref{closed},
and from the vanishing of all the other
pentagram invariants, that
$1+O_{2m}=3O_n=3$.  Hence $O_{2m}=2$.  Likewise
$E_{2m}=2$.
\endproof

\subsection{Proof of Theorem \ref{hyper}}

We start with a degenerate PolyLine.  This
PolyLine determines a labelling of a single
row of 
the triangulation $\tau$, considered in \S 6.
Combining Theorem
\ref{lift} and Lemma \ref{van} we see that
we can find a circulent condensation $C$
which translates into our pentagram labelling.

\begin{center}
\psfig{file=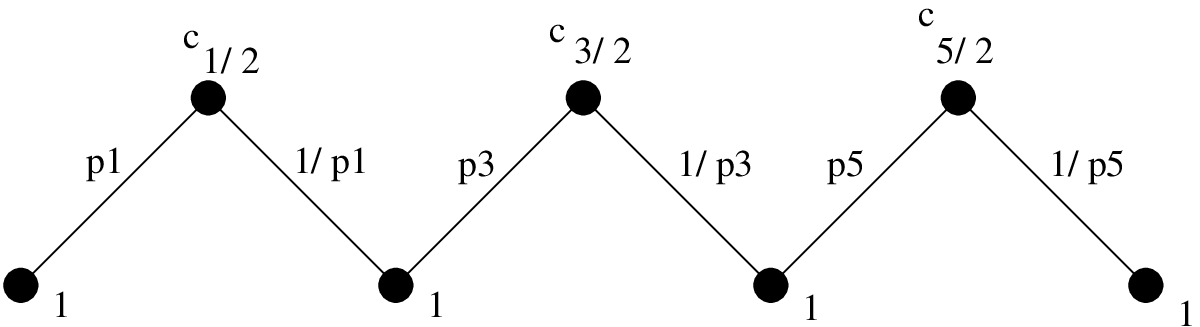}
Figure 5.5
\end{center}

Referring to the proof of Theorem \ref{lift} we
can use the fact that $p_jq_{j+1}=1$ to arrange
so that $r_jr_{j+1}=1$ for all $j$.  This allows
to pick our $c$ labelling so that the row below the
edge labels is identically $1$, as shown in
Figure 5.5.  Figure 5.5 is supposed to be periodic
in the horizontal direction.

We pull our $c$ labelling back to give a labelling of
$T$, the tiling of $\R^3$ by  octahedra.
The bottom row of dots pulls back to an infinite
horizontal plane of dots labelled by $1$.  The next layer up
pulls back to give the labelling of the next
horizontal plane, as shown in Figure 5.6.
This labelling is periodic with respect to
horizontal and vertical translation by $n$.
Each $n \times n$ block is a circulent matrix.
For generic choice of variables, this
circulent matrix will have nonzero determinant.

\begin{center}
\psfig{file=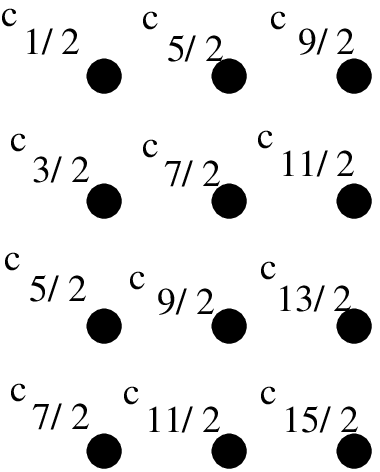}
Figure 5.6
\end{center}

When we develop the condensation upwards, we
arrive precisely at a sandwich condensation of
width $n$.   Going back to the planar picture,
we see that the $n$th horizontal row of our
vertex labelling is a constant labelling,
as shown in Figure 5.7.

\begin{center}
\psfig{file=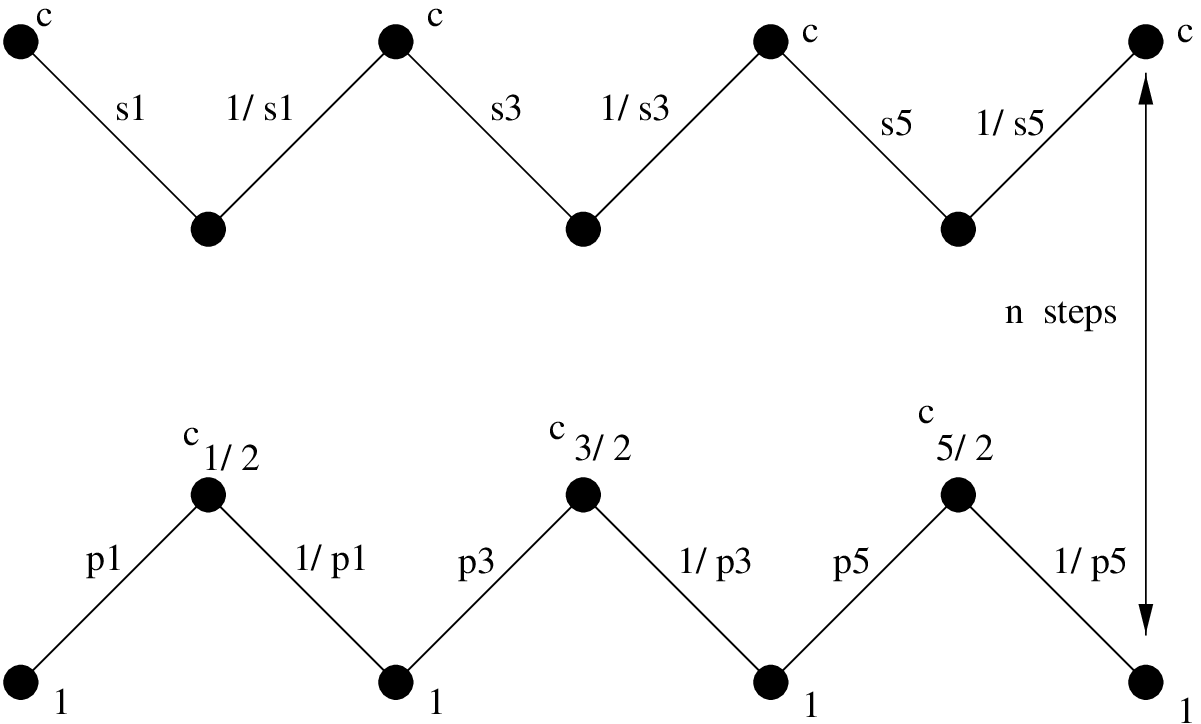}
Figure 5.7
\end{center}

Since the top later of vertices is all constant, the
layer of edges directly below it must have labels
$s_1,t_2,s_3,t_4...$ where $t_2=1/s_1$, $t_4=1/s_3$, etc. 
By Lemma \ref{recip} the PolyPoint corresponding
to this row is degenerate.  That is, the vertices
of this PolyPoint lie on a pair of lines.
Translating this information back into the
language of the pentagram map, we see that
it corresponds exactly to the statement
of Theorem \ref{hyper}.

\newpage

\section{Proof of Theorem \ref{precise}}
\label{independence}

\subsection{Proof modulo the Vanishing Lemma}

In this chapter we prove that the
pentagram invariants,
constructed in \S 2.1, are algebraically independent.
We take $n$ odd, so that there are $n+1$
dihedral invariants.  The even case is very similar.

Given a polynomial map 
$f: \C^{2n} \to \C$, let
$\nabla f=(\partial f/\partial z_1,...,\partial f/\partial z_{2n})$.

\begin{lemma}
\label{C}
Suppose there is a $p \in \C^{2n}$ such
that $\nabla f_1(p),...,\nabla f_{n+1}(p)$ are
linearly independent over $\C$.  Then
$f_1,...,f_{n+1}$ are algebraically independent.
\end{lemma}

\startproof
For the sake of contradiction assume that there is some
nontrivial $F \in \Z[x_1,...,x_{n+1}]$ such
that $F(f_1,...,f_{n+1})=0$.  
By continuity, there is an open set $U$ such
that $\nabla f_1(q),...,\nabla f_{n+1}(q)$ are
linearly independent for any point $q \in U$.
Since $F$ is a nontrivial polynomial map there
is some $q \in U$ such that $dF(q) \not = 0$.
By the chain rule, $\nabla f_k \in {\rm ker\/}(dF)$
for all $k$.
These vectors span $\C^{n+1}$, forcing
$dF(q) =0$, a contradiction.
\endproof

Let $\omega=\exp(2 \pi i/ n)$. For $v \in [1,(n-3)/2]$
let $\Lambda_v$ be the collection of all sequences
$\{s_1,...,s_v\} \subset \{2,3,...,n-2\}$ such that 
$s_j \leq s_{j+1}+2$ for all $j$.   
Define
\begin{equation} \label{exp}
\lambda_v=\sum_{I \in \Lambda_v}\omega^I; \hskip 40 pt
\omega^I=\prod_{j=1}^v \omega^{s_j};
\hskip 15 pt I=(s_i,...,s_v).
\end{equation}
In the next section we will prove
that $\lambda_v \not = v$ for all
$v \in [1,(n-3)/2]$.  In this section
we take this result, which we call
the Vanishing Lemma, for granted.

We have two embeddings of $\C^n$ into $\C^{2n}$.
Namely \begin{equation}
o(z_1,...,z_n)=(z_1,0,z_2,0,...,z_n,0); \hskip 15 pt
e(z_1,...,z_n)=(0,z_1,0,z_2,...,0,z_{n}).
\end{equation}
Define $H_k=O_k \circ o=E_k \circ e$.
The function $H_n$ is defined in $\C^n$.

\begin{lemma}
\label{alg2}
$\nabla H_1(p),...,\nabla H_{(n-1)/2}(p), \nabla H_n(p)$ are
linearly independent at
$p=(\omega,\omega^2,...,\omega^n)$.
\end{lemma}

\startproof
Since $H_{i+1}$ is cyclically invariant, and homogeneous
of degree $i+1$, the gradient $\nabla H_i$ is homogeneous
of degree $i$.  Furthermore, the $(j+1)$th entry of
$\nabla H_i$ is obtained from the $j$th entry by shifting
the indices of the variables by $1$.  These two
facts imply that
$\nabla H_{i+1}(p)=\mu_i V_i,$ where
$$\mu_i=\frac{\partial H_{i}}{\partial z_n}\bigg|_p; \hskip 30 pt
V_i=(\omega^i,\omega^{2i},...,\omega^{ni}).$$
The vectors $V_1,...,V_{\frac{n-1}{2}},V_n$ are
certainly linearly independent over $\C$.
It suffices to prove all the
$\mu_i$ are nonzero.

As $H_{1}=z_1+...+z_n$ and
$H_{n}=z_1...z_n$, we have $\mu_1=\mu_n=1$.
For the intermediate values of $i$, the terms in $H_{i}$ have the form
$(-1)^i z_{k_1},...,z_{k_i}$, where successive or
repeating indices are not allowed.  Note, in particular,
that $...z_1z_n...$ never occurs, because the notion of
succession is reckoned cyclically.  The
terms in $\partial H_i/\partial z_n$ therefore have the form
$$(-1)^iz_{k_1},...,z_{k_{i-1}}; \hskip 30 pt
k_{\alpha} \leq k_{\alpha+1}+2 \hskip 10 pt \forall \alpha.$$
Also, the terms $z_{n-1}$, $z_n$ and $z_1$ do not occur.
Thus, we see that $\mu_i=\pm \lambda_{i-1}$.  Since
$i \leq (n-1)/2$, the Vanishing Lemma says that $\mu_i \not = 0$.
\endproof

\begin{lemma}
\label{limm}
Let $p_t=(t\omega,t\omega,t\omega^2,t\omega^2,...
       t\omega^n,t\omega^n)$.  Then
$$\lim_{t \to 0} t^{1-k} \nabla O_k(p_t)=
o(\nabla H_k(p)); \hskip 20 pt
\lim_{t \to 0} t^{1-k} \nabla E_k(p_t)=
e(\nabla H_k(p)).$$
\end{lemma}

\startproof
We will derive the first equation, the second being similar.
Given any point $q=(x_1,...,x_{2n})$ define 
\begin{equation}
\label{hom}
S_t(q)=(t^{-1}x_1,tx_2,t^{-1}x_3,tx_4,t^{-1}x_5,...,tx_{2n}).
\end{equation}
We have \begin{equation}
\label{ind1}
\lim_{t \to 0} S_{t}(p_t)=
(\omega,0,\omega^2,0,...,\omega^n,0)=o(p).
\end{equation}
Let $\partial_jO_k$ be the $j$th partial derivative
of $O_k$.
For any point $q$ we have the general
homogeneity: \begin{equation}
O_k(S_t(q))=t^{-k}(O_k(q)). \end{equation}
Hence
\begin{equation}
\label{ind2}
\partial_jO_k(S_{t}p_t)=t^{-k_j}\partial_jO_k(p_t)
\hskip 20 pt k_j=k+(-1)^j.
\end{equation}
Combining Equations \ref{ind1} and \ref{ind2} we have
\begin{equation}
\label{ind4}
\lim_{t \to 0} t^{1-k} \partial_jO_k(p_t)=
\left\{\matrix{\partial_jO_k(o(p))&&&{\rm j\ odd\/} \cr \cr 0
&&&{\rm j\ even\/}}\right.
\end{equation}
Finally, \begin{equation}
\label{ind5}
\partial_{2j+1}(O_k(o(p))=o(\partial_j H_k(p)).
\end{equation}
Equations \ref{ind4} and \ref{ind5} together establish
the first equation.
\endproof

We claim that the vectors
$\nabla O_1(p_t),...,\nabla E_n(p_t)$ are
independent for some $t$.  Otherwise there are
functions $a_{k,t}$ and $b_{k,t}$ such that
\begin{equation}
\label{ind7}
\sum a_{k,t} \nabla O_k(p_t)+
\sum b_{k,t} \nabla E_k(p_t)=0;
\hskip 30 pt
\max(|a_{k,t}|,...,|b_{n,t}|)=1.
\end{equation}
We let $t \to 0$.
Taking a subsequence, we can arrange that
$\lim a_{k,t}=a_{k,0}$ and
$\lim_{b,t}=b_{k,0}$ for all relevant $k$, with
at least one limit being nonzero.

We multiply Equation \ref{ind7} by $t^{1-k}$
and take the limit using Lemma \ref{limm} to
obtain
\begin{equation}
\sum a_{k,0}\ o(\nabla H_k(p))+
\sum b_{k,0}\ e(\nabla H_k(p))=0.
\end{equation}
The subspaces $o(\C^n)$ and
$e(\C^n)$ are orthogonal and the
vectors $\{\nabla H_k(p)\}$ are
linearly independent.  This is a contradiction.
Hence $\nabla O_1(p_t),...,\nabla E_n(p_t)$ are
algebraically independent for some value of $t$.
Lemma \ref{C} now completes our algebraic
independence proof, modulo the Vanishing Lemma.

\subsection{Proof of The Vanishing Lemma}
\label{independence proof}

We begin with some algebraic preliminaries.
Say that an {\it adapted measure\/} is a 
positive measure $\tau$, with integer sized
atoms, supported in the $n$th roots of unity.
(Here $n$ is fixed, as above.)
Each adapted measure $\tau$ is encoded by a
non-decreasing sequence 
$I=\{s_1,...,s_k\}$.  The $j$th root
of unity has $\tau$-mass $m$ iff $j$ appears
$m$ times in $I$.
We define $\langle \tau \rangle=\omega^I$, as in Equation \ref{exp}.
By convention, the $\langle \emptyset \rangle=1$.
We define the product
$\tau_1 \cdot \tau_2$ to be the measure
obtained by adding $\tau_1$ and $\tau_2$ together.
(We will reserve the $+$ symbol
for another purpose.)
Note that 
$\langle \tau_1 \cdot \tau_2 \rangle = \langle \tau_1
\rangle \langle \tau_2 \rangle$.  
For $m \in \N$ we define $m \tau=
\tau \cdot ... \cdot \tau$, a total of $m$ times.

Let $M$ denote the free abelian group
generated by the adapted measures.  A
typical element of $M$ is a finite formal
sum $\sigma=\sum m_i \tau _i$.
We define the {\it evaluation map\/}
 \begin{equation}
\langle \sigma \rangle=
\sum \langle \tau_i \rangle^{m_i}.
\end{equation}  We make $M$ into a ring
using the product rule
\begin{equation}
\big(\sum_i m_i \sigma_i\big) 
\big(\sum_j n_j \sigma_j\big)=
\sum_{i,j} m_in_j( \sigma_i \cdot \tau_j)
\end{equation}
The ring $M$ is the {\it group ring\/} generated
by the adapted measures.  The evaluation
map is a ring homomorphism from $M$ to $\C$.

If $A \subset S^1$ is an arc and $v \in \N$ is
an integer, let $\Psi(A,v) \subset M$ denote the
sum, taken over all adapted measures
which have mass $v$ and are supported in $A$.
Using the notation from the Vanishing Lemma,
let $A_v \subset S^1$ be the open arc,
containing $-1$, whose endpoints are
$\omega^v$ and $\omega^{-v}$.

\begin{lemma}[Centrally Symmetric Compression]
$\lambda_v=\langle\Psi(A_v,v)\rangle.$ 
\end{lemma}

\startproof
Let $\Lambda_v$ be the set of sequences used to
define $\lambda_v$ in Equation \ref{exp}.
If $I=(s_1,...,s_v) \in \Lambda_v$, then let
$\phi(I)$ be the summand of $\Psi(A_v,v)$ indexed by
$(s_1+v-1,s_2+v-3,s_2+v-5,...,s_{v-2}-v+5,s_{v-1}-v+3,s_v-v+1).$
This map is a bijection which
preserves the total sum of the elements in $I$,
so that $\omega^I=\langle \phi(I) \rangle$. 
This lemma now follows from the definitions
of the two relevant quantities.
\endproof

Our proof breaks down into two main cases, which we
treat in turn.

\subsubsection{Case 1: $v<n/4$}

We say that a measure is {\it sparse\/} if it
assigns at most mass $1$ to any given point.
For any pair $(A,m)$ let $\Psi'(A,m)$ denote the
formal sum of sparse mass-$m$ measures, supported in
$A$. Let $A^c=S^1-A$.

\begin{lemma} [Binomial Theorem]
$\langle\Psi(A_v,v)\rangle=
\langle\Psi'(A_v^c,v)\rangle$.
\end{lemma}

\startproof
We write $A=A_v$, $A^c=A^c_v$,
$\Psi=\Psi(A_v,v)$, and
$\Psi'=\Psi'(A^c_v,v)$.  
Let $\Phi_j$ be the
formal sum of all mass $v$ measures whose
support intersects $A^c$ in exactly $j$ points.
Note that $\Phi_0=\Psi$ and $\Phi_v=\Psi'$.
Let $\Theta_j$ denote the formal sum of
all mass $j$ adapted measures.  By symmetry
$\langle \Theta_j \rangle=0$ for $j>0$. 
Let $\Delta_j=\Psi_j(A^c,j)$ be the formal
sum of sparse adapted measures of mass $j$ which
are supported in $A^c$.  Note
that $\Delta_v=\Psi'$. 

Suppose that $k \in \{0,....,v-1\}$.  If
$j \geq k$ and $\tau$ is a summand of
$\Phi_j$ there are exactly $j$ choose $k$
ways to write 
$\tau=\tau_1 \cdot \tau_2$, where
$\tau_1 \in \Delta_k$ and
$\tau_2 \in \Theta_{v-k}.$
The point is that we can choose the support of
$\tau_1$ to be any $k$-element subset of
the $A^c$-support of $\tau$. 
This way of counting things gives the relation:
\begin{equation}
\Delta_k \Theta_{v-k}=\sum_{j=k}^v 
\left(\begin{array}{c} j \\k \end{array} \right) \Phi_j,
\end{equation}
for $k=0,...,v-1$.
Combining the previous equation with
a familiar corollary of the binomial theorem,
\begin{equation}
\sum_{k=0}^{v-1} (-1)^k \Delta_k \Theta_{v-k}=
\Phi_0 + (-1)^v \Phi_v.
\end{equation}
Since
$\langle \Delta_k \Theta_{v-k} \rangle=
\langle \Delta_k \rangle \langle \Theta_{v-k} \rangle=
0.$ we have
$\langle \Psi \rangle=
\langle \Phi_0 \rangle= \pm
\langle \Phi_v \rangle= 
\langle \Psi' \rangle$.
\endproof

Since $v<n/4$ we have
$\Re(z)>0$ for all $z \in A_v$.
We will use induction 
to show that $\langle \Psi'(A_v,w) \rangle>0$ for all $v,w \geq 1$.
Let $\underline \omega^v$ be the mass $1$ measure
supported on $\omega^v$.  If $\tau$ is a
mass $w$ sparse measure supported in $A_v^c$ then
the support of $\tau$ intersects
$\{\omega^v,\omega^{-v}\}$ in $0$, $1$, or $2$
points.  Thus
\begin{equation}
\label{indu}
\Psi'(A_v^c,w)=\left\{\matrix{\Psi'(A^c_{v-1},w) \cr + \cr
(\underline \omega^v+\underline \omega^{-v}) 
\cdot \Psi'(A^c_{v-1},w-1) \cr + \cr
(\underline \omega^v \cdot \underline \omega^{-v}) \cdot 
\Psi'(A^c_{v-1},w-2).}\right\}
\end{equation}
At least one term on the right is nontrivial.  From
\begin{equation}
\label{indu2}
\langle \underline \omega^v+\underline \omega^{-v}\rangle=
2 \Re(\omega^v)>0; \hskip 30 pt
\langle \underline \omega^v \cdot \underline \omega^{-v} \rangle=1.
\end{equation} 
and induction, any nontrivial term on the
right hand side of Equation \ref{indu} evaluates
to a positive number.   Therefore, the
left hand side evaluates to a positive
number as well.

\subsubsection{Case 2: $v \geq n/4$}

For each integer $w \in (0,n/4]$ we choose an open arc $B_w$, 
invariant under complex conjugation, such that
$-1 \in B_w$ and
there are exactly $w$ $n$th roots of
unity contained in $B_w$.
Let $\Psi(w,k',k)$ denote the formal sum of
adapted mass $k$ measures $\mu$ such that
$\mu$ is supported in $B_w$ and
$\mu(B_w-B_{w-2}) \leq k'$.

Our goal is to show that
$\langle \Psi(w,v,v) \rangle \not = 0$,
where $w$ is the number of $n$th roots of
unity in $A_v$.   
We order the triples $(w,k',k)$ lexicographically.
We will show inductively that
$\langle \Psi(w,k',k) \rangle>0$ if $k$ is even
and $\langle \Psi(w,k',k)<0$ if $k$ is odd.
(These sums are real, by symmetry.)

If $k=1$ then 
$\langle \Psi(w,k',k) \rangle$ is the sum of numbers
all of which have negative real part, so that
$\langle \Psi(w,k',k) \rangle<0$ in this case.
Also,
$\langle\Psi(1,k,k) \rangle=(-1)^k.$
Henceforth we assume that $w \geq 2$ and
$k \geq 2$.   Since $w \geq 2$ there are
two $n$th roots of unity
$\alpha_1$ and $\alpha_2=\overline \alpha_1$
in $B_w-B_{w-2}$.

Suppose $w=2$.
A simple counting argument gives
$$\Psi(w,k,k)=(\underline \alpha_1 + \underline \alpha_2) \cdot
\Psi(v,k-1,k-1) + \underline \alpha_1 \cdot \underline \alpha_2
\cdot \Psi(v,k-2,k-2).$$
Note that $\alpha_1+\alpha_2<0$.
By induction, both terms on the right have the
desired sign when evaluated.
Henceforth we assume that $w \geq 3$.

Suppose that $k'=1$.
A counting argument gives
$$\Psi(w,1,k)=\Psi(w-2,k,k)+
(\underline \alpha_1+\underline \alpha_2) \cdot
\Psi(w-2,k-1,k-1)$$
Again, we note that $\alpha_1+\alpha_2<0$.
Since
$w \geq 3$ both terms on the right have the 
desired sign when evaluated.

Suppose that $k'=2$.
A counting argument gives
$$\Psi(w,2,k)=\Psi(w-2,k,k)+
(\underline \alpha_1+\underline \alpha_2) \cdot \Psi(1,k-1)+
\underline \alpha_1 \cdot \underline \alpha_2
\cdot \Psi(w-2,k-2,k-2).$$
By induction, all terms on the right have
the desired sign when evaluated.

Suppose that $k' \geq 3$.
A counting argument gives
$$
\label{induct}
\Psi(w,k',k)= \left\{ \matrix{
\Psi(w-2,k'-2,k) \cr + \cr
(\underline \alpha_1 + \underline \alpha_2) \cdot
\Psi(w,k'-1,k-1) \cr + \cr
\underline \alpha_1 \cdot \underline \alpha_2
\cdot \Psi(w,k'-2,k'-2)} \right\} .
$$
By induction, all three terms on the right
have the desired sign when evaluated.

This completes our proof.

\newpage

\section{References}

[{\bf RR\/}] D.P. Robbins and H. Rumsey, \newline
{\it Determinants and Alternating Sign Matrices\/} \newline
Advances in Mathematics {\bf 62\/} (1986) 
\newline
\newline
[{\bf S1\/}] R. Schwartz, {\it The Pentagram Map\/} \newline
Journal of Experimental Mathematics {\bf 1\/} (1992)  pp. 85-90
\newline
\newline
[{\bf S2\/}] R. Schwartz, {\it Recurrence of the Pentagram Map\/}, \newline
Journal of Experimental Mathematics (2001)
\newline
\newline
[{\bf S3\/}] R. Schwartz, \newline
 {\it Desargues Theorem, Dynamics, and Hyperplane Arrangements\/}  \newline
Geometriae Dedicata (2001)
\newline
\newline
[{\bf W\/}] S. Wolfram,
{\it The Mathematica Book\/},
Fourth Edition,
Cambridge University Press (1999)

\newpage

\end{document}